\documentclass[11pt]{article}
\usepackage[utf8]{inputenc}
\usepackage{float,cite,bm,graphicx,amsmath,amsthm,amssymb,subcaption,caption,bbold,amstext,array,authblk,mathtools,xcolor}
\usepackage[small]{titlesec}
\usepackage{epstopdf}
\usepackage[nottoc]{tocbibind}
\newcolumntype{L}{>{$}l<{$}} 
\usepackage[margin=0.8in]{geometry}
\usepackage{setspace}
\mathtoolsset{showonlyrefs}

\newcommand{\bx}{\bm{x}}

\newcommand{\bz}{\bm{z}}

\newcommand{\bg}{\bm{g}}

\newcommand{\bff}{\bm{f}}
\newcommand{\by}{\bm{y}}

\newcommand{\dd}{\mathrm{d}}
\newcommand{\R}{\mathbb{R}}

\newcommand{\p}{\partial}

\newtheorem{theorem}{Theorem}
\newtheorem{remark}{Remark}

\theoremstyle{definition}

\begin{document}
	\title{Geometric integration of ODEs using multiple quadratic auxiliary variables}
	\author{Benjamin K Tapley\thanks{email: bentapley@hotmail.com}}
	\affil{Department of Mathematical Sciences, The Norwegian University of Science and Technology, 7491 Trondheim, Norway}
	\date{\today}
	\maketitle
	\begin{abstract}
		We present a novel numerical method for solving ODEs while preserving polynomial first integrals. The method is based on introducing multiple quadratic auxiliary variables to reformulate the ODE as an equivalent but higher-dimensional ODE with only quadratic integrals to which the midpoint rule is applied. The quadratic auxiliary variables can subsequently be eliminated yielding a midpoint-like method on the original phase space. The resulting method is shown to be a novel discrete gradient method. Furthermore, the averaged vector field method can be obtained as a special case of the proposed method. The method can be extended to higher-order through composition and is illustrated through a number of numerical examples. 
	\end{abstract}
	\section{Introduction}
	In this paper, we are concerned with the numerical solution of autonomous ODEs in $n$ dimensions
	\begin{equation}
		\dot{\bx} = \bff(\bx), \quad \bx \in \R^n, \label{ODE}
	\end{equation}
	that possess functions $H_i(\bx):\R^n\rightarrow\R$ for $i=1,...,k<n$ satisfying $\dot{H}_i(\bx) = \bff(\bx)\cdot\nabla H_i(\bx)=0$, where the dot denotes $\frac{\mathrm{d}}{\mathrm{d}t}$. If such functions exist, they are called \textit{first integrals} or \textit{invariants} of \eqref{ODE}. In the $k=1$ integral case, one can write \eqref{ODE} in skew-gradient form $\bff(\bx)=S(\bx)\nabla H(\bx)$ for some skew-symmetric matrix $S(\bx)=-S(\bx)^T$ \cite{mclachlan1998unified}. When $S(\bx) = {\tiny \left(\begin{array}{cc} 0 & -I \\I & 0 \\ \end{array}\right)}$ for identity matrix  $I$, then \eqref{ODE} is called Hamiltonian. The fact that first integrals are preserved along the exact solution of \eqref{ODE} is a reflection of the underlying physical laws that govern qualitative features of the system's dynamics, therefore the preservation of these integrals by numerical methods is important, especially for long-term simulations \cite{hairer2006geometric}.  Numerical methods that inherent physical properties of the systems they model are referred to as \textit{geometric methods} and, in the case at hand, the physical property is integral-preservation (often referred to as energy-preservation). Hamiltonian systems are a large and important class of conservative ODE systems that govern many physical phenomena. Consequently, many well designed energy-preserving methods have been proposed for solving such systems, for example, Hamiltonian boundary value methods \cite{brugnano2015reprint} and Runge-Kutta type methods \cite{miyatake2014energy,ranocha2020relaxation,iserles2000preserving}. Another prevalent class of energy-preserving methods are the discrete gradient methods \cite{gonzalez1996time,mclachlan1999geometric}. The idea behind the discrete gradient method is to find a function $\bar{\nabla}H(\bx',\bx):\R^{n}\times\R^{n}\rightarrow\R^{n}$ satisfying
	\begin{equation}\label{dg}
		\bar{\nabla}H(\bx',\bx)\cdot(\bx'-\bx) = H(\bx')-H(\bx)\\
	\end{equation}
	along with the consistency condition $\bar{\nabla}H(\bx,\bx) = {\nabla}H(\bx)$. Such a function is a discrete approximation to $\nabla H(\bx)$. An integral-preserving discrete gradient method can then be constructed by solving
	\begin{equation}
		\frac{\bx'-\bx}{h} = \tilde{S}\left(\bx',\bx,h\right)\bar{\nabla}H(\bx',\bx),
	\end{equation}
	where $\bx'$ is an approximation to the exact solution $\bx(h)$ from initial conditions $\bx$ and $\tilde{S}\left(\bx,\bx,0\right)={S}\left(\bx\right)$. Due to the above properties, it can be shown that the discrete gradient method is energy-preserving, i.e., $H(\bx')=H(\bx)$. Some of the first and most widely used discrete gradients are attributed to Itoh and Abe \cite{itoh1988hamiltonian} and Gonzalez \cite{gonzalez1996time} and we refer also to \cite{mclachlan1999geometric} for a good exposition on discrete gradients and their application in geometric integration. Another prevailing example of a discrete gradient method is the averaged vector field (AVF) method \cite{quispel2008new}, which when applied to canonical Hamiltonian systems with constant $S$ yield the following energy-preserving method
	\begin{equation}
		\frac{\bx'-\bx}{h} = S\!\int_{0}^{1}\nabla H((1-\xi)\bx + \xi \bx')\dd \xi.
	\end{equation}
	 Due to its simplicity, the AVF method has been the building block of many interesting methods for canonical Hamiltonian ODEs \cite{hairer2010energy, celledoni2009energy}. However, such methods are not only limited to ODEs. The study of energy-preserving numerical methods for ODEs often serves as inspiration for many energy-preserving methods for conservative \textit{partial} differential equations (PDEs) including discrete variational methods \cite{matsuo2001dissipative,furihata2019discrete}, AVF based methods \cite{celledoni2012preserving,gong2014some}, continuous Runge-Kutta type methods \cite{li2015general} and fast linearly implicit methods based on polarization \cite{dahlby2011general} or discrete gradients \cite{eidnes2019linearly}. Other efficient and well documented energy-preserving method for conservative PDEs are the scalar auxiliary variable approach \cite{shen2018scalar,kemmochi2021scalar} and the invariant energy quadratization approach \cite{yang2017numerical}. \\
	
	
	Here, we present a novel numerical method that preserves polynomial integrals of ODEs. The method is straight-forward to implement and can be summarised in the following steps:
	\begin{enumerate}
		\item Introduce quadratic auxiliary variables to reformulate the polynomial integrals as quadratic functions.
		\item Reformulate the ODE in terms of the auxiliary variables.
		\item Discretise the reformulated ODE with the midpoint rule, which preserves the quadratic integrals.
		\item Eliminate the auxiliary variables from the discretisation to arrive at a midpoint-like method on the original phase space. 
	\end{enumerate} 
	The method can be applied to systems that possess multiple first integrals, but it's efficiency is limited by the computation of a $k+1$-skew-symmetric tensor. We report that a similar method was developed concurrently to the present paper named the quadratic auxiliary (QAV) method  \cite{gong2021new}. Here, the authors propose a very efficient method for the KdV PDE based on introducing a single quadratic auxiliary variable to reduce the cubic Hamiltonian into a quadratic, then integrating the system with a symplectic Runge-Kutta method. The numerical results reported suggest that the method, together with a well designed structure-preserving implementation lead to a very efficient algorithm for the KdV equation. In contrast, the present paper describes how to take advantage of \textit{multiple} quadratic auxiliary variables to address the time integration of conservative ODEs with arbitrary polynomial first integrals. We therefore refer to the proposed method as the multiple quadratic auxiliary variable (MQAV) method. \\
	
	As the MQAV method is essentially the midpoint rule applied to the reformulated ODE it is therefore implicit, symmetric, A-stable and roughly the same computational cost of an implicit Runge-Kutta method of equal order. The order of the method can be increased using symmetric composition or by diagonally implicit symmetric Runge-Kutta methods. In general, the MQAV method is not uniquely defined and one can generate a family of MQAV methods that depend on a number of free parameters that can be tuned to suit application. Furthermore, the MQAV method fits within the framework of existing geometric numerical integration methods. In particular, we will show that the MQAV method is a new class of discrete gradient. Moreover, when a certain choice of free parameters are made, the MQAV method surprisingly produces the AVF method as a special case. Lastly, we add that while the focus is on time integration of ODEs, the method extends straightforwardly to semi-discretised PDEs in a similar fashion to the methods proposed in \cite{gong2021new, celledoni2012preserving}.\\
		
	The paper is organised as follows. In section \ref{sec:Quartic Ham} we introduce the MQAV method for the quartic integral case. Here, we present and give theoretical results on the MQAV method then illustrate it using the planar quartic oscillator as an example. Extending the method to preserve arbitrary polynomial integrals follows straightforwardly and is presented in sections \ref{sec:polynomial integrals} and \ref{sec:multiple integrals}, respectively. In section \ref{sec:higher order} we show how to develop higher-order MQAV methods. Section \ref{sec:numerical examples} presents numerical examples of a Hamiltonian system, a Nambu system and the Toda lattice. Concluding remarks are given in section \ref{conclusion}.

	%
	\section{Integral-preserving numerical integration using multiple quadratic auxiliary variables}
	\subsection{ODEs with a quartic first integral}\label{sec:Quartic Ham}
	Let $\R_d[\bx]$ be the class of polynomials of degree $d$ in the variables $\bx\in\R^{n}$ and consider an ODE of the form \eqref{ODE} with a first integral $H(\bx)$
	\begin{equation}
		\dot{\bx} = S(\bx)\nabla H(\bx), \label{SGODE}
	\end{equation}
	where $S(\bx)=-S(\bx)^T$ is a skew-symmetric matrix.  \\ 
	
	We start by introducing the quadratic auxiliary variables 
	\begin{equation}
		y_{i,j} = x_ix_j \quad \text{for}\quad j\le i, \quad i=1,...,n
	\end{equation} 
	and define a reduced-degree Hamiltonian $\tilde{H}(\bx,\by)\in\R_2[\bx,\by]$ that is at most \textit{quadratic} in the new variables $\bx$ and $\by$ and satisfies the consistency condition \begin{equation}
		\tilde{H}(\bx,\by(\bx))=H(\bx).
	\end{equation} 
	Here, with a slight abuse of notation, we have let $\by = (y_{1,1},y_{1,2},...,y_{n,n-1},y_{n,n})^T \in \R^{m}$ denote the $m=n(n+1)/2$ independent quadratic auxiliary variables and $\by(\bx) = (x_1x_1,x_1x_2,...,x_nx_{n-1},x_nx_n)^T$ (i.e., with the argument $\bx$ written explicitly) denote the variables $\by$ as a function of the original phase space variables $\bx$. Using the chain rule, the ODE in these new variables now read
	\begin{equation}
		\begin{aligned}\label{xlp}
			\dot{\bx} =& S(\bx) \nabla \tilde{H}(\bx,\by)=S(\bx) \left(\frac{\p}{\p \bx} \tilde{H}(\bx,\by) +  \frac{\p \by(\bx)}{\p \bx}^T \frac{\p }{\p \by}\tilde{H}(\bx,\by)\right):=\bff(\bx,\by)\\
			\dot{\by} =& \frac{\p \by(\bx)}{\p \bx}\bff(\bx,\by):=\bg(\bx,\by)
		\end{aligned}
	\end{equation}
	We will refer to $\bff(\bx,\by)$ as the \textit{MQAV formulation} of $\bff(\bx)$ with respect to $\tilde{H}(\bx,\by)$, which also satisfies the consistency condition $\bff(\bx,\by(\bx)) = \bff(\bx)$. Furthermore, we will refer to the collective ODE system \eqref{xlp} as the \textit{extended} ODE. We now make some observations. 
	\begin{remark}\label{rem:1}
		The extended ODE \eqref{xlp}:
		\begin{enumerate}
			\item[(a)] is of the following skew gradient form for $\bz = (\bx^T,\by^T)^T$  
			\begin{equation}
				\dot{\bz} = K(\bz)\nabla\tilde{H}(\bz), \quad K(\bz) = -K(\bz)^T= \left(\begin{array}{cc}
					S(\bx) & -\left(\frac{\p \by(\bx)}{\p \bx}S(\bx)\right)^T\\
					\frac{\p \by(\bx)}{\p \bx}S(\bx) & \frac{\p \by(\bx)}{\p \bx}S(\bx)\frac{\p \by(\bx)}{\p \bx}^T\end{array}\right).\\
			\end{equation}
			\item[(b)] has $\tilde{H}(\bx,\by)$ as first integral. This follows directly from (a) due to the skew-symmetry of $K(\bz)$. 
			\item[(c)] possesses the $m$ additional quadratic integrals
			\begin{equation}
				H_{i,j}(\bx,\by) = x_ix_j - y_{i,j}, \quad i=1,...,n,\quad j\le i.
			\end{equation}
			This can be easily seen by differentiating with respect to $t$
			\begin{equation}
				\dot{H}_{i,j} = \dot{x}_ix_j + \dot{x}_jx_i - \dot{y}_{i,j} = 0.
			\end{equation}
			\item[(d)] is solved by the solution of \eqref{SGODE} with the consistent initial conditions $y_{i,j}\big|_{t=0} = x_ix_j\big|_{t=0}$.
		\end{enumerate} 
	\end{remark}
	Due to remark \ref{rem:1}d, an integral-preserving method of the extended ODE \eqref{xlp} is equivalent to an integral-preserving method on the original ODE \eqref{SGODE} as long as the initial conditions $y_{i,j}\big|_{t=0} = x_ix_j\big|_{t=0}$ are satisfied.
	Moreover, the $m$ induced quadratic  integrals $H_{i,j}(\bx,\by)$ together with the reduced-degree Hamiltonian $\tilde{H}(\bx,\by)$ are quadratic, therefore any symplectic Runge-Kutta method (i.e., one with vanishing stability matrix \cite{cooper1987stability}) will preserve all such integrals. However, numerically integrating the extended ODE \eqref{xlp} presents a significant drawback, namely that the phase space dimension has increased from $n$ to $n+m=n(n+3)/2$. Due to the fact that $m$ scales quadratically with $n$ this quickly becomes computationally prohibitive for large $n$ especially due to the fact all symplectic Runge-Kutta methods are implicit. However, this issue can be circumvented as follows. First we implement the (symplectic) midpoint rule given by
	\begin{align}
		\frac{\bx' - \bx}{h} = \bff\left(\frac{\bx' + \bx}{2},\frac{\by' + \by}{2}\right), \label{1}\\
		\frac{\by' - \by}{h} = \bg\left(\frac{\bx' + \bx}{2},\frac{\by' + \by}{2}\right), \label{2}
	\end{align}
	where $\bx'\approx \bx(h)$ denotes the solution of the midpoint rule from initial conditions $\bx=\bx(0)$ and similarly for $\by$. Due to the fact that the midpoint rule preserves all quadratic invariants, we have preservation of the induced quadratic  integrals $H_{i,j}(\bx',\by') = H_{i,j}(\bx,\by)$, which form a set of $m$ linear equations for $\by'$ meaning its solution can be directly found in terms of $\bx'$ and therefore avoiding the need to numerically solve \eqref{2}. That is, by substituting $y'_{i,j} = x'_ix'_j$ into \eqref{1} we are left with an $n$ dimensional implicit equation for $\bx'$ alone. This is given by
	\begin{equation}\label{rdmp}
		\frac{\bx' - \bx}{h} = \bff\left(\frac{\bx' + \bx}{2},\frac{\by(\bx') + \by(\bx)}{2}\right).
	\end{equation}
	Doing so circumvents the issue of the increased phase space dimension, thus yielding an efficient method that preserves any quartic Hamiltonian. The method \eqref{rdmp} will be referred to as the \textit{MQAV midpoint method}. 
	\begin{remark}
		We observe that the MQAV midpoint method:
		\begin{enumerate}
			\item[(a)] is integral-preserving, that is, $H(\bx')=H(\bx)$.
			\item[(b)] is order two.
			\item[(c)] is symmetric, that is, it is invariant under the transformation $\bx\leftrightarrow\bx'$ and $h\rightarrow -h$.
			\item[(d)] is affinely equivariant with respect to the variables $(\bx,\by)$.
			\item[(e)] possesses a B-series with respect to the extended ODE vector field. 
			\item[(f)] is A-stable.
		\end{enumerate} 
	\end{remark}
	These observations follow from the fact that the MQAV midpoint rule is the conventional midpoint rule applied to the extended ODE and therefore inherits similar properties but on the extended vector field. We also note that while the midpoint rule is symplectic, the MQAV midpoint rule is not due to the fact that reformulating the ODE into its extended form destroys the symplectic structure of the original ODE in the Hamiltonian case. Furthermore, an advantage of having a B-series expansion is that the method automatically preserves affine symmetries and linear integrals of the ODE. However in this setting, these symmetries and integrals of the original ODE need also to be reflected by the extended ODE to be preserved by the MQAV method. \\ 

		A key ingredient of the MQAV method is to first reduce the degree of an integral that is quartic in $\bx$ to one that is quadratic in $\bx$ and $\by$. It turns out that there is not a unique way to represent an integral in its reduced-degree form. Consider an arbitrary quartic first integral in the variables $\bx\in\R^n$
		\begin{equation}\label{ham4}
			H(\bx) = \sum_{i} \alpha_i x_{a_i}x_{b_i}x_{c_i}x_{d_i},
		\end{equation}	
		where $a_i,b_i,c_i,d_i$ can take values from 0 to $n$ and $x_0:=1$ to include monomials of degree less than four. Then the most general way to reduce the Hamiltonian to quadratic is
		\begin{equation}\label{rdh4}
			\tilde{H}(\bx,\by) = \sum_{i} \alpha_i \left(\beta^{[i]}_1 y_{a_i,b_i}y_{c_i,d_i}+\beta^{[i]}_2 y_{d_i,a_i}y_{b_i,c_i}+\beta^{[i]}_3 y_{a_i,c_i}y_{b_i,d_i}\right),
		\end{equation}		
		where $\beta^{[i]}_1,\beta^{[i]}_2$ and $\beta^{[i]}_3$ are free parameters satisfying $\beta^{[i]}_1+\beta^{[i]}_2+\beta^{[i]}_3=1$ for consistency. This lets us develop a family of methods where the parameters $\beta^{[i]}_j$ can be tuned to improve accuracy or stability, for example. With this general reduced-degree representation of $H(\bx)$, one can show how this fits into the framework of the discrete gradient methods. \\
		
		\begin{theorem}\label{thmdg}
			The MQAV gradient corresponding to \eqref{rdh4} is a discrete gradient. That is, 
			\begin{equation}
				\bar{\nabla}_{MQAV}\tilde{H}(\bx',\bx) = \nabla\tilde{H}\left(\frac{\bx'+\bx}{2},\frac{\by(\bx')+\by(\bx)}{2}\right)
			\end{equation}
			satisfies equation \eqref{dg} when $\beta^{[i]}_1+\beta^{[i]}_2+\beta^{[i]}_3=1$ and is consistent. 
		\end{theorem}
		\begin{proof}
			Let  
			\begin{equation}
				T_{a_i,b_i,c_i,d_i;j}(\bx,\by) := \frac{\partial}{\partial x_j} y_{a_i,b_i}y_{c_i,d_i}\\ 
				= \left(\delta_{j,a_i}x_{b_i}+\delta_{j,b_i}x_{a_i}\right)y_{c_i,d_i} + y_{a_i,b_i}\left(\delta_{j,c_i}x_{d_i}+\delta_{j,d_i}x_{c_i}\right)
			\end{equation}	
			for Dirac-delta $\delta_{i,j}$. Letting $\bar{\bx} = \frac{\bx'+\bx}{2}$ and $\bar{\by} = \frac{\by(\bx')+\by(\bx)}{2}$ we can write the MQAV discrete gradient as
			\begin{equation}
				\nabla_j\tilde{H}(\bar{\bx},\bar{\by}) = \sum_{i} \alpha_i \left(\beta^{[i]}_1T_{a_i,b_i,c_i,d_i;j}(\bar{\bx},\bar{\by})+\beta^{[i]}_2T_{d_i,a_i,b_i,c_i;j}(\bar{\bx},\bar{\by})+\beta^{[i]}_3T_{a_i,c_i,b_i,d_i;j}(\bar{\bx},\bar{\by}) \right).
			\end{equation}
			Contracting the above with $x'_j-x_j$ and simplifying the right hand side gives
			\begin{equation}
				\sum_{j}\nabla_j\tilde{H}(\bar{\bx},\bar{\by})(x'_j-x_j) = \sum_{i} \alpha_i \left(\beta^{[i]}_{1}+\beta^{[i]}_{2}+\beta^{[i]}_{3}\right) \left(x'_{a_{i}} x'_{b_{i}} x'_{c_{i}} x'_{d_{i}} - x_{a_{i}} x_{b_{i}} x_{c_{i}} x_{d_{i}} \right)=H(\bx')-H(\bx)\\
			\end{equation}
			which, together with the consistency condition, completes the proof. 
		\end{proof}
	To our knowledge, it appears that $\bar{\nabla}_{MQAV}\tilde{H}(\bx',\bx)$ is a novel type of discrete gradient. We now recall another popular discrete gradient corresponding to the AVF method
	\begin{equation}\label{AVF}
		\bar{\nabla}_{AVF}H(\bx',\bx) = \int_0^1 \nabla H((1-\xi)\bx+\xi\bx') \dd\xi.
	\end{equation}	
	It turns out that the AVF discrete gradient is a special case of the MQAV midpoint method when applied to Hamiltonian ODEs with quartic Hamiltonian due to the following. \\
	
		\begin{theorem}\label{thm:avf}
			For quartic first integrals, the MQAV discrete gradient corresponding to \eqref{rdh4} is the AVF discrete gradient when $\beta^{[i]}_1=\beta^{[i]}_2=\beta^{[i]}_3=1/3$. That is, 
			\begin{equation}
				\bar{\nabla}_{MQAV}\tilde{H}(\bx',\bx)=\bar{\nabla}_{AVF}H(\bx',\bx) 
			\end{equation}
		\end{theorem}
		\begin{proof}
			This can be seen by direct computation, similarly to that of Theorem \ref{thmdg}. 
		\end{proof}
	
		While is it obvious that the MQAV method is a Runge-Kutta method on the extended space (i.e., the midpoint rule) it is less obvious whether or not it yields a Runge-Kutta method on the original space. However in \cite{celledoni2014minimal} it is shown that the AVF method coincides with the Runge-Kutta method whose Butcher coefficients correspond to the nodes and weights of a Gaussian quadrature applied to the integral in \eqref{AVF}. Therefore, due to \ref{thm:avf} we can infer that the MQAV discrete gradient also coincides with a Runge-Kutta method. An interesting question is whether or not this holds for all choice of the free parameters $\beta_{i}^{[j]}$. While this section considers quartic integrals, we have seen that the MQAV method can also reproduce the AVF method as a special case and is a discrete gradient method for octic integrals, which is evidence to support the claims that in the general $H(\bx)\in\R_d[\bx]$ case: (1) the AVF method arises as a special case of the MQAV method and (2) the MQAV gradient is a discrete gradient, although we omit these results. \\
		 
		
	\subsubsection{Example: The quartic oscillator}
\begin{figure}
		\centering
		\includegraphics[width=0.5\linewidth]{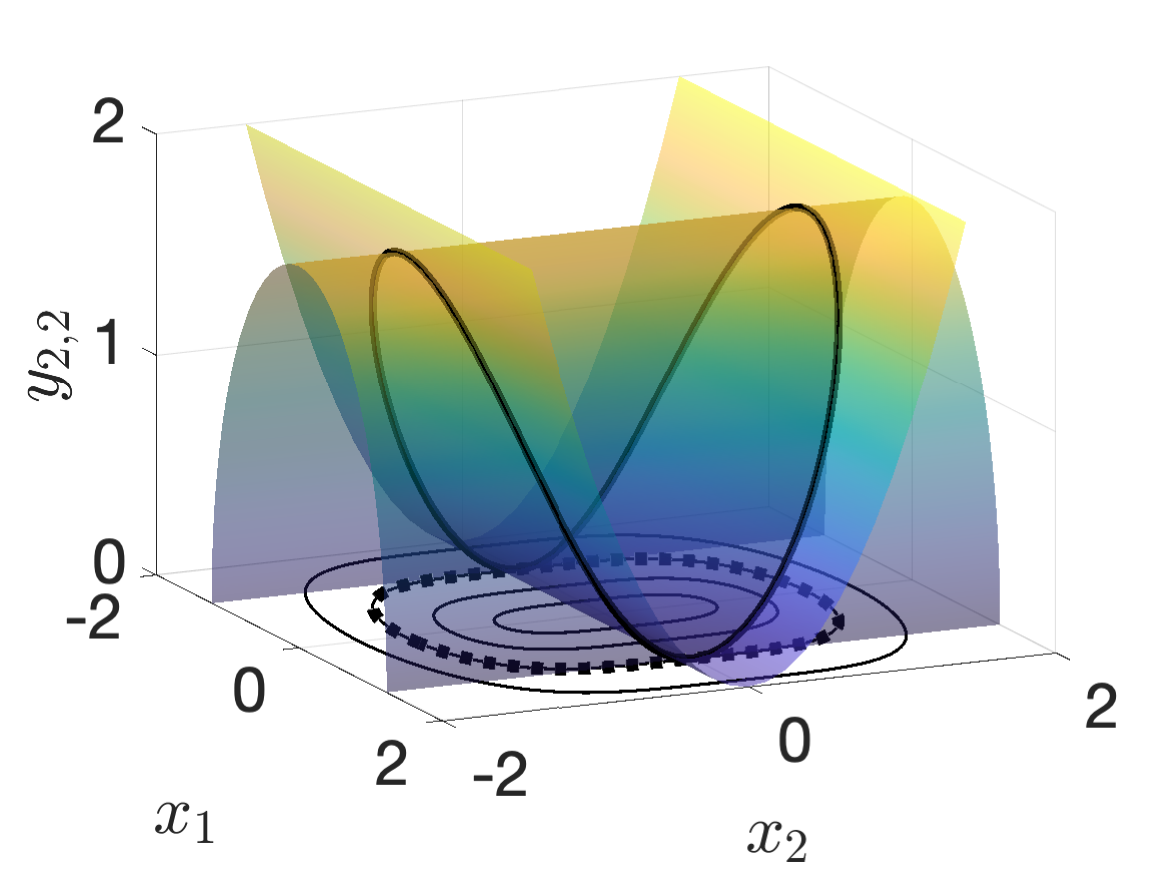}
		\caption{The numerical solution of the midpoint method applied to the 3 dimensional extended ODEs (thick black line), the MQAV midpoint method applied to the original 2 dimensional Hamiltonian ODE (dotted black line), the level sets of the original quartic Hamiltonian ${H}(x_1,x_2)$ (thin black lines) and the iso-surfaces of the reduced-degree Hamiltonian $\tilde{H}(x_1,x_2,y_{2,2})$ and the induced quadratic  integral $H_{2,2}(x_1,x_2,y_{2,2})$ (colored surfaces). The simulation uses the initial conditions $x_1(0)=1$, $x_2(0)=1$ and $y_{2,2}(0)=x_2(0)^2$. Note that the MQAV midpoint method is the projection of the midpoint method applied to the three dimensional extended ODEs onto the $x_1-x_2$ plane, which coincides with the level set of ${H}(x_1,x_2)$. }
		\label{fig:phasequadrics}
	\end{figure}

	Consider the quartic oscillator with Hamiltonian 
	\begin{equation}
		H(x_1,x_2) = \frac{1}{2}x_1^2 + \frac{1}{4}x_2^4\in\R_4[x_1,x_2]
	\end{equation} 
	that corresponds to the ODE
	\begin{align}
		\dot{x}_1 = & -x_2^3,\\
		\dot{x}_2 = & x_1.
	\end{align}
	Now introduce the variable $y_{2,2} = x_2^2$ and define the following reduced-degree Hamiltonian 
	\begin{equation}
		\tilde{H}(x_1,x_2,y_{2,2}) = \frac{1}{2}x_1^2 + \frac{1}{4}y_{2,2}^2\in\R_2[x_1,x_2,y_{2,2}].
	\end{equation} 
	Note that this satisfies the consistency condition $\tilde{H}(x_1,x_2,x_{2}^2) = {H}(x_1,x_2)$. The corresponding extended ODE system for $x_1$, $x_2$ and $y_{2,2}$ (i.e., equation \eqref{xlp}) is therefore 
	\begin{equation}
		\begin{aligned}\label{h1}
			\dot{x}_1 = & -x_2y_{2,2},\\
			\dot{x}_2 = & x_1,\\
			\dot{y}_{2,2} = & 2x_1x_2,\quad\text{where}\quad y_{2,2}\big|_{t=0}=x_2^2\big|_{t=0},
		\end{aligned}
	\end{equation}
	to which we apply the midpoint rule (i.e., equations \eqref{1} and \eqref{2})
	\begin{equation}
		\begin{aligned}\label{mp1}
			\frac{x_1'-x_1}{h} = & -\left(\frac{x_2'+x_2}{2}\right)\left(\frac{y_{2,2}'+y_{2,2}}{2}\right),\\
			\frac{x_2'-x_2}{h} = & \left(\frac{x_1'+x_1}{2}\right),\\
			\frac{y_{2,2}'-y_{2,2}}{h} = & 2\left(\frac{x_1'+x_1}{2}\right)\left(\frac{x_2'+x_2}{2}\right).
		\end{aligned}
	\end{equation}
	The ODE \eqref{h1} possesses the integrals $\tilde{H}(x_1,x_2,y_{2,2})$ and $H_{2,2}(x_1,x_2,y_{2,2}) = y_{2,2} - x_2^2$, both of which are quadratic and are therefore preserved by the midpoint rule. This means $y_{2,2}' - x_2'^2 = y_{2,2} - x_2^2$ and due to the fact that the initial conditions satisfy $y_{2,2}=x_2^2$, we have $y_{2,2}' = x_2'^2$. Substituting this into equation \eqref{mp1} we derive the MQAV midpoint rule (i.e., equation \eqref{rdmp})
	\begin{equation}
		\begin{aligned}\label{1rdmp1}
			\frac{x_1'-x_1}{h} = & -\left(\frac{x_2'+x_2}{2}\right)\left(\frac{x_2'^2+x_2^2}{2}\right),\\
			\frac{x_2'-x_2}{h} = & \left(\frac{x_1'+x_1}{2}\right),
		\end{aligned}
	\end{equation}
	which is equivalent to the midpoint rule applied to the extended system \eqref{mp1} projected onto the $\bx$ coordinates. This concept is depicted in figure \ref{fig:phasequadrics}. Here, the numerical solution using the midpoint rule applied to the three dimensional extended ODE \eqref{mp1} is given by the thick black line while the MQAV midpoint rule \eqref{1rdmp1} is given by the dashed line and is confined to the $x_1-x_2$ plane. The isosurfaces of the reduced-degree Hamiltonian $\tilde{H}(x_1,x_2,y_{2,2})$ and the induced integral $H_{2,2}(x_1,x_2,y_{2,2})$ are also displayed as coloured surfaces. We see that the numerical flow of the midpoint rule \eqref{mp1} is confined to the intersection of these two quadrics. Moreover, the numerical flow of the MQAV midpoint method is none other than the projection of this curve onto the $x_1-x_2$ plane. This projection also coincides with the level set of the original quartic Hamiltonian $H(x_1,x_2)$. Lastly, we add that in this case there is a unique representation of the reduced-degree Hamiltonian, therefore this method is identical to the AVF method. 
	
	\subsection{Preservation of polynomial integrals}\label{sec:polynomial integrals}
	 We now focus on the case where $H(\bx)\in\R_d[\bx]$ is a polynomial integral of degree $d$ and show how to extend the MQAV midpoint rule to preserve such an integral.\\
	
	Similarly to the previous section, we start by defining a quadratic reduced-degree integral $$\tilde{H}(\bx,\by^{[1]},...,\by^{[n_y]}) \in \R_2[\bx,\by^{[1]},...,\by^{[n_y]}]$$ for $n_y\ge \log_2(d)-1$ by introducing the change of variables inductively defined by $y^{[0]}_{i}=x_i$ and $y^{[k]}_{\mathbf{i},\mathbf{j}}=y^{[k-1]}_{\mathbf{i}}y^{[k-1]}_{\mathbf{j}}$, where $\mathbf{i}$ is a multi-index. For example, $y^{[1]}_{i,j} = x_ix_j$ and $y^{[2]}_{i,j,k,l} = y^{[1]}_{i,j}y^{[1]}_{k,l}$, where the latter can be abbreviated using the multi-index notation by $y^{[2]}_{\mathbf{i},\mathbf{j}} = y^{[1]}_{\mathbf{i}}y^{[1]}_{\mathbf{j}}$. Note that the degree of $y^{[k]}_{\mathbf{i}}(\bx)$ is $2^{k}$ in the variables $\bx$. Then, by the chain rule, one can write the gradient operator of the reduced-degree integral
	\begin{align}
		\nabla\tilde{H}(\bx&,\by^{[1]},...,\by^{[n_y]}) \\ &=\left(\frac{\p}{\p \bx}  +  \left(\frac{\p \by^{[1]}(\bx)}{\p \bx}\right)^T \frac{\p }{\p \by^{[1]}}+...+ \left(\frac{\p \by^{[n_y]}(\bx)}{\p \bx}\right)^T \frac{\p }{\p \by^{[n_y]}}\right) \tilde{H}(\bx,\by^{[1]},...,\by^{[n_y]}) \\ &= \sum_{k=0}^{n_y} \left(\frac{\p \by^{[k]}(\bx)}{\p \bx}\right)^T \frac{\p }{\p \by^{[k]}}\tilde{H}(\bx,\by^{[1]},...,\by^{[n_y]})\label{total deriv}
	\end{align}
	where $\frac{\p \by^{[k]}(\bx)}{\p \bx}$ is the Jacobian derivative matrix of $\by^{[k]}(\bx)$ with respect to $\bx$, also computed by the chain rule. 
	The MQAV ODE for $\bx$ with respect to $\tilde{H}(\bx,\by^{[1]},...,\by^{[n_y]})$ is therefore
	\begin{align}\label{ODE1}
		\dot{\bx} = S(\bx)\nabla \tilde{H}(\bx,\by^{[1]},...,\by^{[n_y]}) :=\bff(\bx,\by^{[1]},...,\by^{[n_y]}),
	\end{align}
	which satisfies the consistency condition $\bff(\bx,\by^{[1]}(\bx),...,\by^{[n_y]}(\bx)) = \bff(\bx)$. Although not needed for our method, the corresponding ODEs for the $\by^{[k]}$ variables that make up the extended system are
	\begin{equation}
		\dot{\by}^{[k]} = \frac{\p \by^{[k]}(\bx)}{\p \bx} \bff(\bx,\by^{[1]},...,\by^{[n_y]})\quad\text{for}\quad k=1,...,n_y.
	\end{equation}
	These, together with $\dot{\bx}$ form an extended system that possesses the induced quadratic  integrals $H_{\mathbf{i},\mathbf{j}} = y^{[k]}_{\mathbf{i},\mathbf{j}} - y^{[k-1]}_{\mathbf{i}}y^{[k-1]}_{\mathbf{j}}$. As each $H_{\mathbf{i},\mathbf{j}}$ as well as the reduced degree Hamiltonian $\tilde{H}$ are all quadratic then the midpoint rule $\bx'$ preserves these integrals. Therefore we can immediately solve for $\left(\by^{[k]}\right)'$ in terms of $\bx'$. Inserting this solution yields the MQAV midpoint rule for polynomial integrals
	\begin{align}\label{rdmp2}
		\frac{\bx'-\bx}{h} =\bff\left(\frac{\bx'+\bx}{2},\frac{\by^{[1]}(\bx')+\by^{[1]}(\bx)}{2},...,\frac{\by^{[n_y]}(\bx')+\by^{[n_y]}(\bx)}{2}\right), 
	\end{align}
	the solution to which satisfies $H(\bx')=H(\bx)$. \\
	
	\subsubsection{Example: the octic oscillator}
	Consider the octic oscillator with Hamiltonian 
	\begin{equation}
		H(x_1,x_2) = \frac{1}{2}x_1^2 + \frac{1}{8}x_2^8\in\R_8[x_1,x_2]
	\end{equation} 
	that corresponds to the ODE
	\begin{equation}
		\begin{aligned}\label{octicHamVF}
			\dot{x}_1 = & -x_2^7,\\
			\dot{x}_2 = & x_1.
		\end{aligned}
	\end{equation}
	Now introduce the variables $y_{2,2}^{[1]} = x_2^2$ and $y_{2,2,2,2}^{[2]} = \left(y_{2,2}^{[1]}\right)^2$ to define the following reduced-degree Hamiltonian 
	\begin{equation}
		\tilde{H}(x_1,x_2,y_{2,2}^{[1]},y_{2,2,2,2}^{[2]}) = \frac{1}{2}x_1^2 + \frac{1}{8}\left(y_{2,2,2,2}^{[2]}\right)^2\in\R_2[x_1,x_2,y_{2,2}^{[1]},y_{2,2,2,2}^{[2]}].
	\end{equation} 
	The MQAV vector field of \eqref{octicHamVF} with respect to $\tilde{H}(x_1,x_2,y_{2,2}^{[1]},y_{2,2,2,2}^{[2]})$ is therefore
	\begin{equation}
		\begin{aligned}\label{octicHamRDVF}
			\dot{x}_1 = & -x_2y_{2,2}^{[1]}y_{2,2,2,2}^{[2]},\\
			\dot{x}_2 = & x_1.
		\end{aligned}
	\end{equation}
	The corresponding extended system possesses the induced quadratic  integrals $H_{2,2} = y_{2,2}^{[1]} - x_2^2$ and $H_{2,2,2,2} = y_{2,2,2,2}^{[2]} - \left(y_{2,2}^{[1]}\right)^2$. Therefore the midpoint method satisfies  $y_{2,2}^{[1]'} = x_2^{'2}$ and $y_{2,2,2,2}^{[2]'} = \left(y_{2,2}^{[1]'}\right)^2 = x_2^{'4}$. Applying the MQAV midpoint rule \eqref{rdmp2} to \eqref{octicHamRDVF} yields 
	\begin{equation}
		\begin{aligned}
			\begin{aligned}\label{octicHamRDMP}
				\frac{x_1' - x_1}{h} = & -\left(\frac{x_2' + x_2}{2}\right)\left(\frac{x_2^{'2} + x_2^2}{2}\right)\left(\frac{x_2^{'4} + x_2^4}{2}\right),\\
				\frac{x_2' - x_2}{h} = & \frac{x_1' + x_1}{2},
			\end{aligned}
		\end{aligned}
	\end{equation}
	the solution to which satisfies $H(x'_1,x_2')=H(x_1,x_2)$. We remark that this method is also identical to the AVF method. 
	

	\subsection{Preservation of multiple integrals}\label{sec:multiple integrals}
	The idea behind the MQAV midpoint rule is simply to introduce auxiliary quadratic variables to write the system as an equivalent extended system that possesses only quadric integrals, then apply the midpoint rule and use the quadratic-integral-preserving property to exactly project the map back onto the original phase space. This can be applied directly to ODEs with multiple first integrals, therefore extending the MQAV midpoint method \eqref{rdmp2} to preserve multiple polynomial integrals is straight forward. This approach yields a similar method to the multiple-integral-preserving discrete gradient methods outlined in \cite{mclachlan1999geometric,norton2015projection}, which we will now briefly outline.\\
	
	 Given an ODE with $k<n$ first integrals $H_i(\bx)\in\R_{d_i}[\bx]$ for $i=1,...,k$, then it is always possible to write the ODE in the following multi-skew-gradient form \cite{quispel1997solving,mclachlan1999geometric}
	\begin{equation}
		\dot{\bx} = S\left(\nabla H_1(\bx) ,\,...\,,\nabla H_k(\bx) \right)
	\end{equation}
	where $S$ is some skew-symmetric $k+1$ tensor depending on $\bx$. That is, $S_{i_0,...,i_a,...,i_b,...i_{k}}=-S_{i_0,...,i_b,...,i_a,...i_{m}}$, for any $a,b=0,...,k$. In general, there are infinitely many solutions for $S$. To find a suitable solution, one can use a symbolic algebra software package such as Maple, use a more heuristic approach outlined in \cite{quispel1997solving,mclachlan1999geometric}, for example or refer to the default formula, which can be written down using the exterior algebra \cite{mclachlan2003spatial}
	\begin{equation}
		S(\bx) = \frac{\bff\wedge\nabla H^{[1]}\wedge ... \wedge\nabla H^{[k]}}{\det(\mathcal{H})}
	\end{equation}
	where $\mathcal{H}_{i,j}=\nabla H^{[i]}\cdot\nabla H^{[j]}$. The ${n \choose k}$ unique tensor components of $S$ can be found by computing the following determinants 
	\begin{equation}\label{s}
		S_{i_0,...,i_k}(\bx) = \frac{1}{\det(\mathcal{H})}\cdot
		\det\left(
		\begin{array}{cccc}
			{f}_{i_0}(\bz) & \partial_{i_0}H^{[1]} & \cdots & \partial_{i_0}H^{[k]} \\
			\vdots & \vdots  &  & \vdots \\
			{f}_{i_k}(\bz) & \partial_{i_k}H^{[1]} & \cdots & \partial_{i_k}H^{[k]} \\
		\end{array}
		\right).
	\end{equation}
	To create an integral preserving method, we start by writing down the $k$ reduced-degree integrals $\tilde{H}_i(\bx,\by^{[1]},...,\by^{[n_y]}) \in \R_2[\bx,\by^{[1]},...,\by^{[n_y]}]$. The MQAV ODE is therefore 
	\begin{equation}
		\dot{\bx} = S\left(\nabla \tilde{H}_1 ,\,...\,,\nabla \tilde{H}_m \right) := \bff(\bx,\by^{[1]},...,\by^{[n_y]})
	\end{equation}
	where the vectors $\nabla \tilde{H_i}$ are functions of $(\bx,\by^{[1]},...,\by^{[n_y]})$ and are calculated by equation \eqref{total deriv}. Then the MQAV midpoint rule from equation \eqref{rdmp2} applied to $\bff(\bx,\by^{[1]},...,\by^{[n_y]})$ preserves all integrals $H_i$ for $i=1,...,m$. Examples are given in section \ref{sec:numerical examples}. We remark that the computation of the $k+1$-tensor $S$ can be costly for high dimensional problems with large values of $k$. 
	\subsection{Higher-order methods}\label{sec:higher order}
	After defining the MQAV vector field and the extended ODE system (e.g., equation \eqref{xlp}) one can apply any symplectic Runge-Kutta method and project the solution onto the original $n$ dimensional phase space to yield an integral-preserving method. However, if we implement higher-order symplectic Runge-Kutta methods, it becomes difficult to avoid calculation of the stage values for the $\dot{\by}^{[k]}$ equations and therefore results in a method that requires the numerical solution of a non-linear system in a much higher dimensional phase space than the original ODE, which is slow. As we have seen in section \ref{sec:Quartic Ham}, the midpoint rule is a special because we can avoid solving the $\by^{[k]}$ equations, and their subsequent stage values, by substituting in the solutions $y'_{i,j}=x'_ix'_j$ that solve the superfluous equations. There are other Runge-Kutta methods that share this property, namely the diagonally-implicit symplectic Runge-Kutta (DISRK) methods. Such methods have Butcher table of the form 
	\begin{equation}\label{SDIRK table}
		\begin{array}{c|cccccc}
			c_1 			& \frac{b_1}{2}  & 0& 0 &\dots & 0 \\
			c_2 			& b_1 & \frac{b_2}{2} & 0 & \dots &  0 \\
			\vdots			 & \vdots & \vdots  &\ddots & \ddots & \vdots \\
			c_{s-1}			 & b_1 & b_2 & \ddots & \frac{b_{s-1}}{2} & 0 \\
			c_s 			& b_1 & b_2 &  \dots& b_{s-1} & \frac{b_s}{2} \\
			\hline
			~				&b_1 & b_2 & \dots & b_{s-1} & b_s \\
		\end{array}
	\end{equation}
	that is, with $a_{ii} = b_i/2$, $a_{ij} = b_j$ for $i<j$ and $a_{ij}=0$ otherwise. These methods are by construction symplectic. Moreover, due to theorem 4.4 in \cite[p. 192 ]{hairer2006geometric} a Runge-Kutta method $\Phi_h$ with coefficients given by \eqref{SDIRK table} is equivalent to compositions of the midpoint rule, which we denote by $\Phi^{M}_h$, with sub-steps $b_ih$. That is, 
	\begin{equation}\label{comp}
		\Phi_h = \Phi^{M}_{b_sh}\circ \dots \circ\Phi^{M}_{b_2h} \circ \Phi^{M}_{b_1h}.
	\end{equation}
	This means that one can construct high-order MQAV Runge-Kutta methods without the need to solve the stage values for the $\by^{[k]}$ equations by taking compositions of the MQAV midpoint rule. The advantage of this is that the MQAV method will inherit the same properties as the underlying Runge-Kutta method when constructed with these coefficients (e.g., A-stability).  \\
	
	In addition to the above DISRK methods, due to the fact that the MQAV midpoint rule is symmetric, one can use composition \cite{mclachlan1995numerical} to create methods of higher-order, in a similar way to equation \eqref{comp} \cite{hairer2006geometric}. This is demonstrated in section \ref{sec:numerical examples}. \\

	\section{Numerical examples}\label{sec:numerical examples}
	This section is comprised of three numerical examples. The first is a planar Hamiltonian ODE, where we compare the MQAV midpoint method against some other geometric integrators for Hamiltonian systems. The second example is of a three dimensional Nambu system possessing two integrals of degree four and eight. Here, we implement higher-order MQAV methods of up to degree eight using DISRK methods and symmetric composition and compare them to the standard methods of equal order. The final example is of the famous Toda lattice system from quantum mechanics. Here, we show that we can numerically integrate the system while preserving all integrals. 
	
	\subsection{A planar quartic-Hamiltonian ODE}
	In this example we will consider a planar Hamiltonian ODE given by
	\begin{equation}\label{qodeeg}
		\left.
		\begin{array}{l}
			\dot{x}_1 = -2\,{x_{1}}^2\,x_{2}-4\,{x_{2}}^3\\
			\dot{x}_2 = 2\,x_{1}\,{x_{2}}^2+x_{1}
		\end{array}\right\} := \bff(\bx)
	\end{equation}
	where the Hamiltonian function is $H(\bx) = \frac{1}{2}x_1^2  + x_2^4 + x_1^2x_2^2$. Then defining a reduced-degree Hamiltonian $\tilde{H}(\bx,\by) = x_1^2/2  + y_{2,2}^2 + \alpha y_{1,1}y_{2,2}+(1-\alpha)y_{1,2}^2$, where  $\alpha\in\mathbb{R}$ is a free parameter, we can derive the corresponding MQAV ODE system using equation \eqref{xlp}
	\begin{equation}
		\left.
		\begin{array}{l}
			\dot{x}_1 =  -2\,(\left( 1-\alpha \right) x_{{1}}y_{{1,2}}+ \,\alpha\,x_{{2}}y_{{1,1}})-4\,x_{{2}}y_{{2,2}} 
			\\ 
			\dot{x}_2 = 2\,(\alpha\,x_{{1}}y_{{2,2}}+\left( 1-\alpha \right)x_{{2}} y_{{1,2}})+x_{{1}}

		\end{array}\right\} := \bff(\bx,\by) \label{qlpode}
	\end{equation}
	We now apply the MQAV midpoint rule, defined by equation \eqref{rdmp}, and set $\alpha=0$ which yields
	\begin{align}\label{rdmpqode}
		\frac{x'_1-x_1}{h} = & -2\,\left(\frac{x_{1}'+x_1}{2}\right)\,\left(\frac{x_{1}'x_{2}'+x_1x_{2}}{2}\right)-4\,\left(\frac{x_{2}'+x_2}{2}\right)\,\left(\frac{x_{2}'x_{2}'+x_{2}x_{2}}{2}\right),\\ 
		\frac{x'_2-x_2}{h} = & 2\,\left(\frac{x_{2}'+x_2}{2}\right)\,\left(\frac{x_1'x_{2}'+x_1x_2}{2}\right)+\left(\frac{x_{1}'+x_1}{2}\right).
	\end{align}	
	\begin{figure}[!h]
		\centering
		\begin{subfigure}{0.32\textwidth}
			\includegraphics[width=\linewidth]{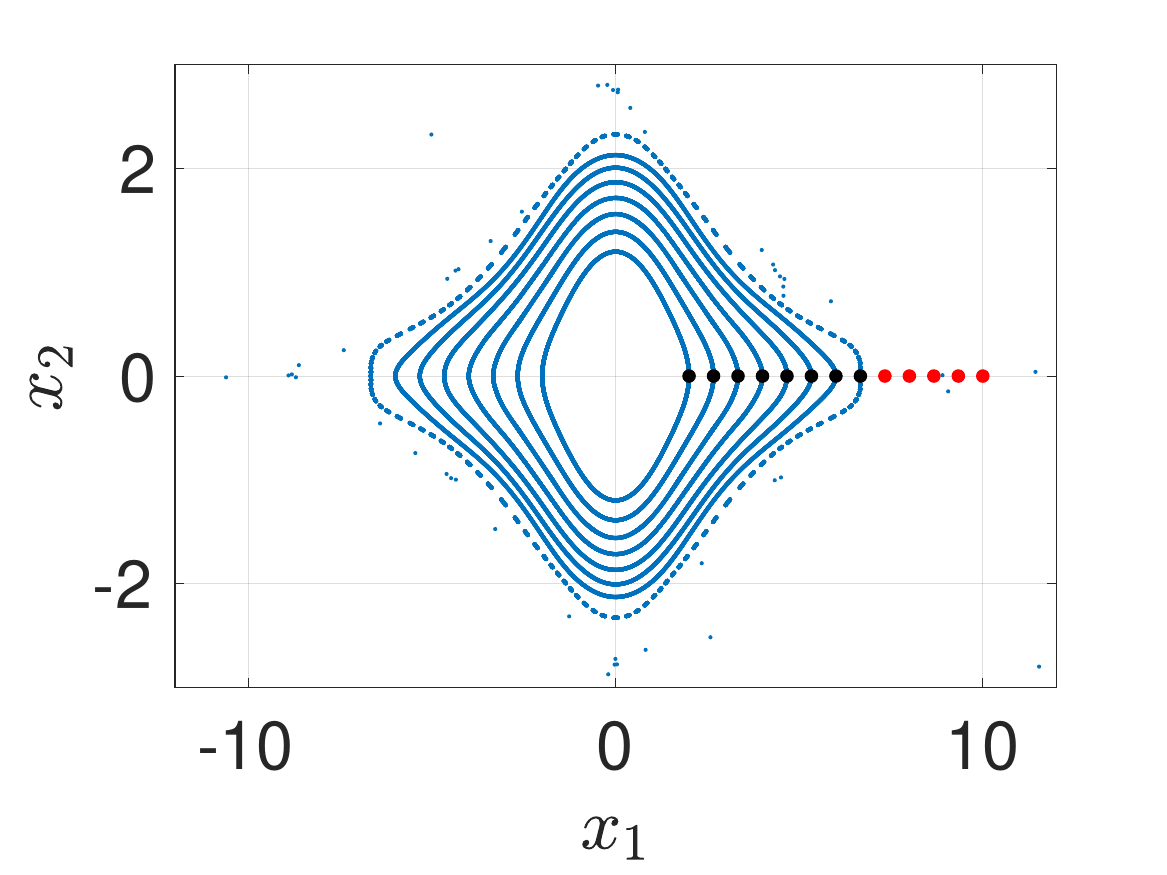}
			\caption{MP2}
			\label{fig:quartic_ham_phaseplot_midpoint_degred}
		\end{subfigure}
		\begin{subfigure}{0.32\textwidth}
			\includegraphics[width=\linewidth]{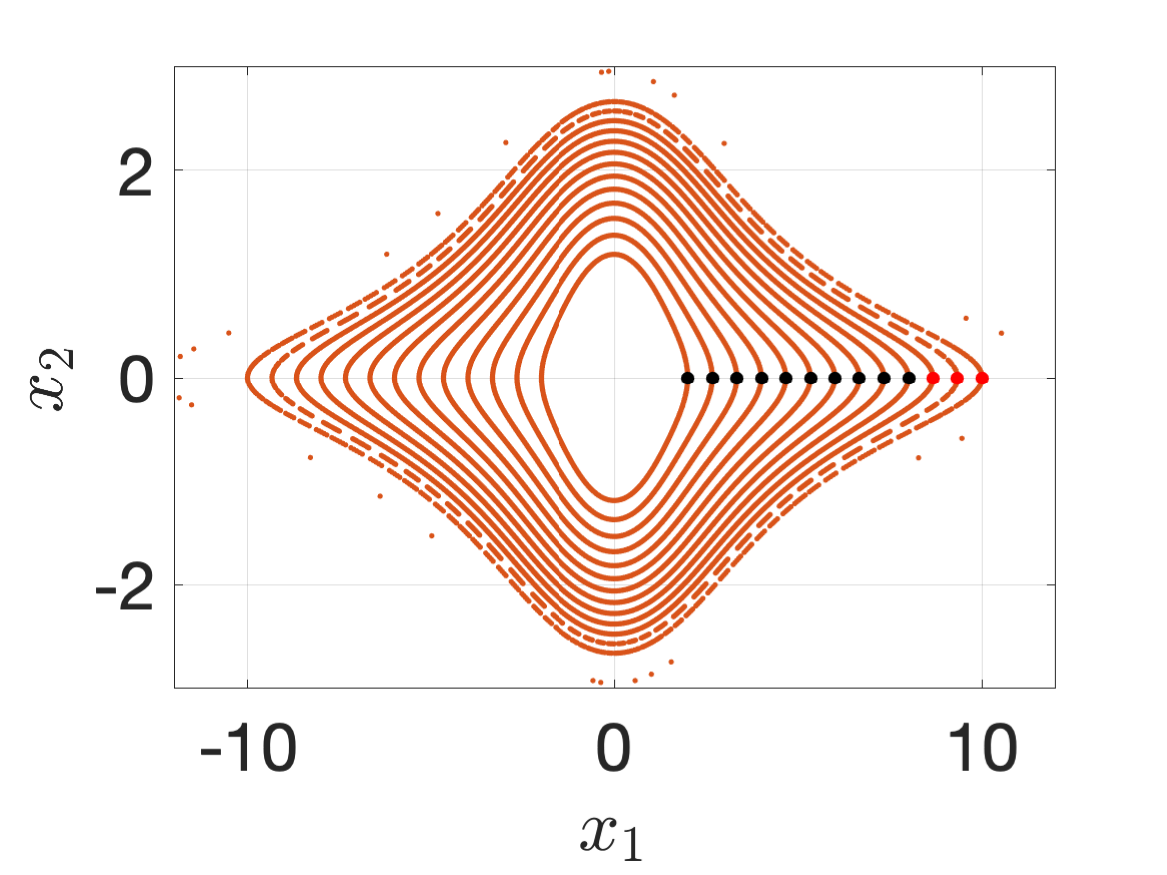}
			\caption{AVF}
			\label{fig:quartic_ham_phaseplot_simpson}
		\end{subfigure}
		\begin{subfigure}{0.32\textwidth}
			\includegraphics[width=\linewidth]{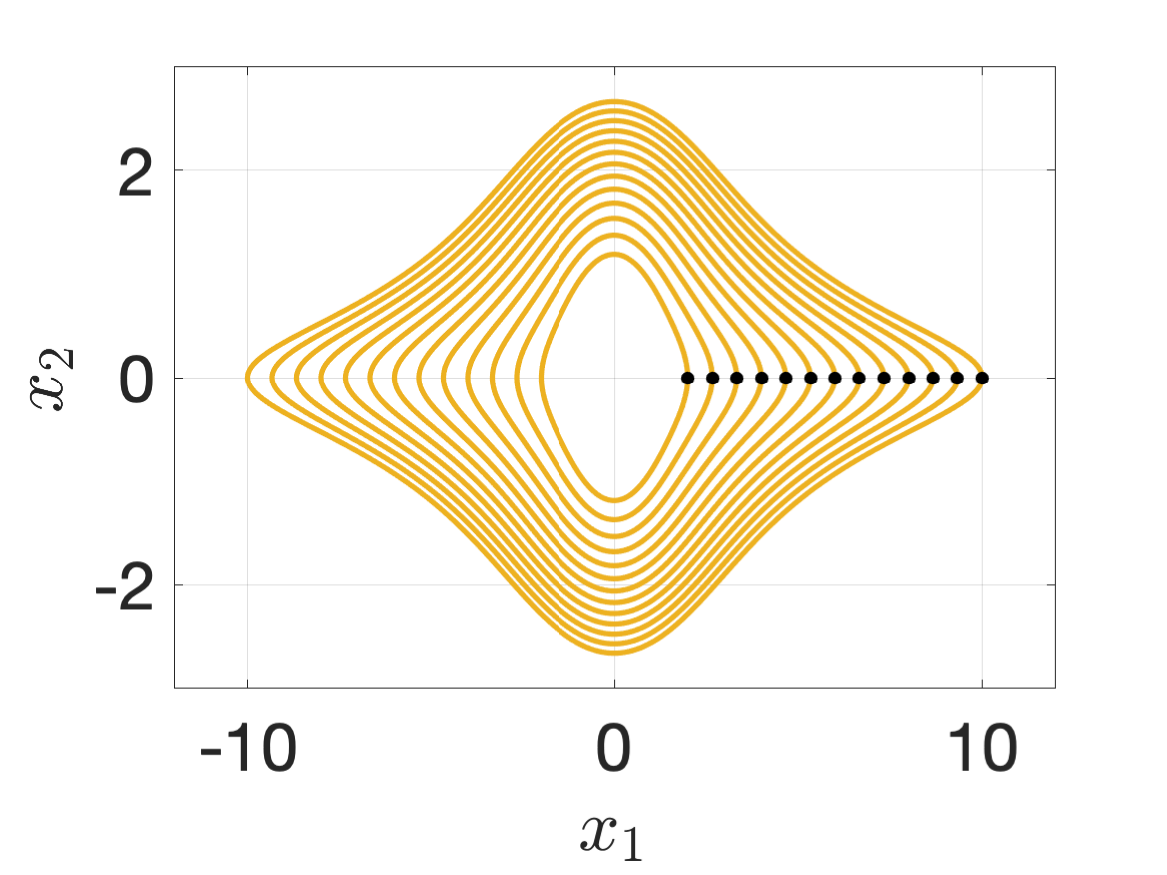}
			\caption{MQAV-MP2}
			\label{fig:quartic_ham_phaseplot_midpoint}
		\end{subfigure}
		
		\begin{subfigure}{0.32\textwidth}
			\includegraphics[width=\linewidth]{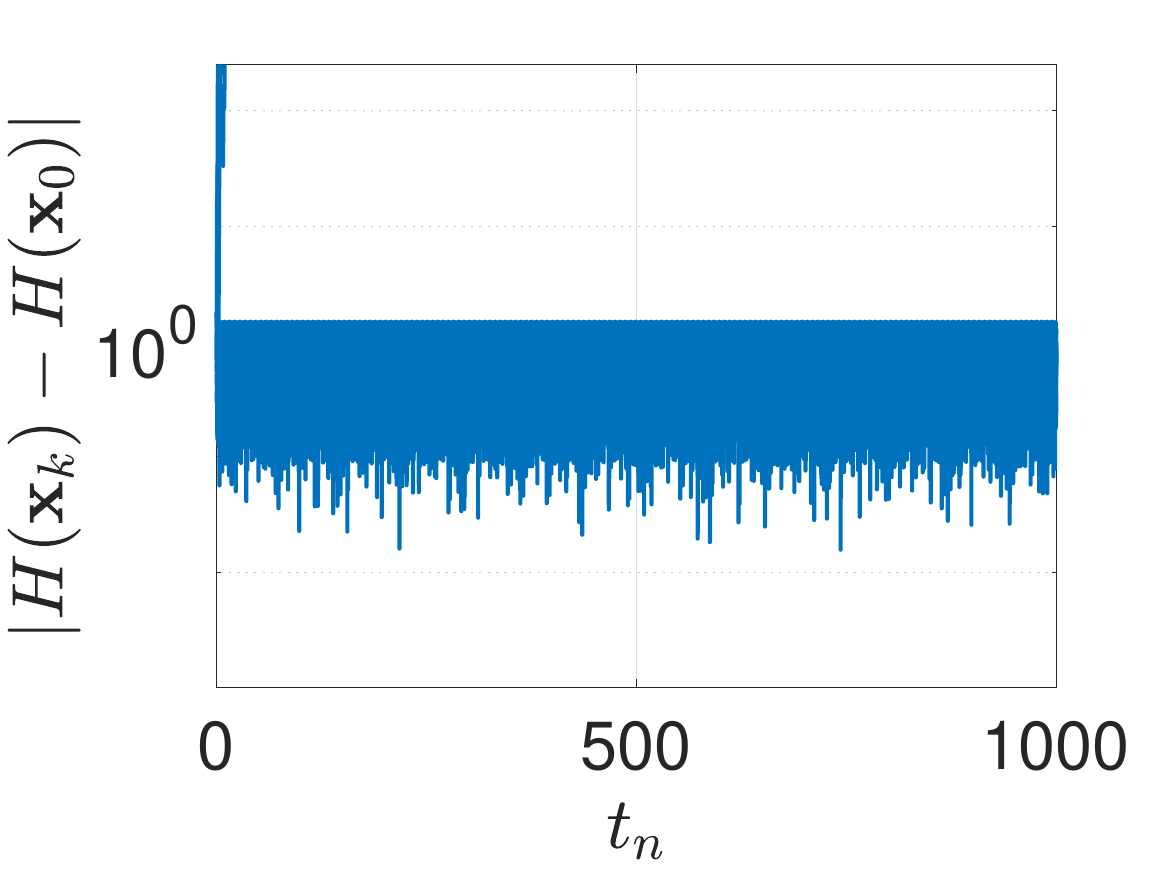}
			\caption{MP2}
			\label{fig:quartic_ham_err1}
		\end{subfigure}
		\begin{subfigure}{0.32\textwidth}
			\includegraphics[width=\linewidth]{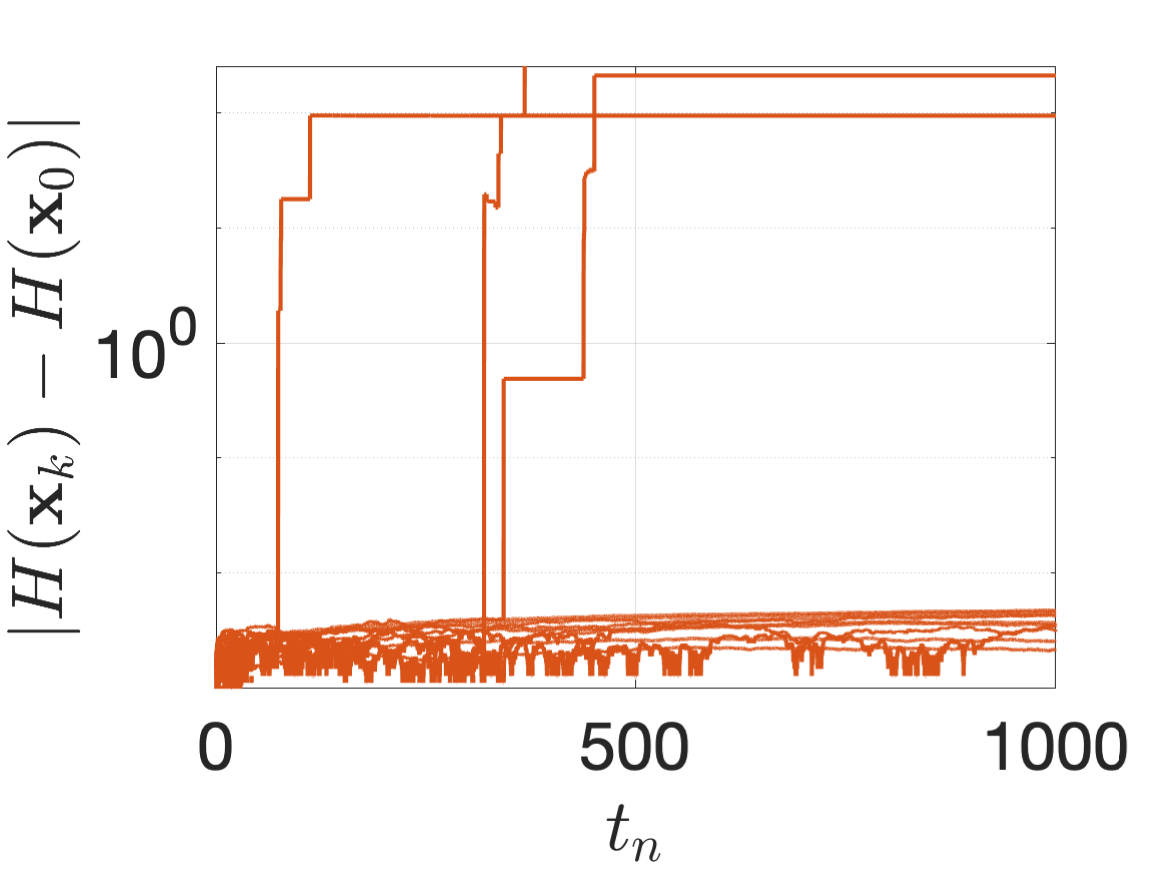}
			\caption{AVF}
			\label{fig:quartic_ham_err2}
		\end{subfigure}
		\begin{subfigure}{0.32\textwidth}
			\includegraphics[width=\linewidth]{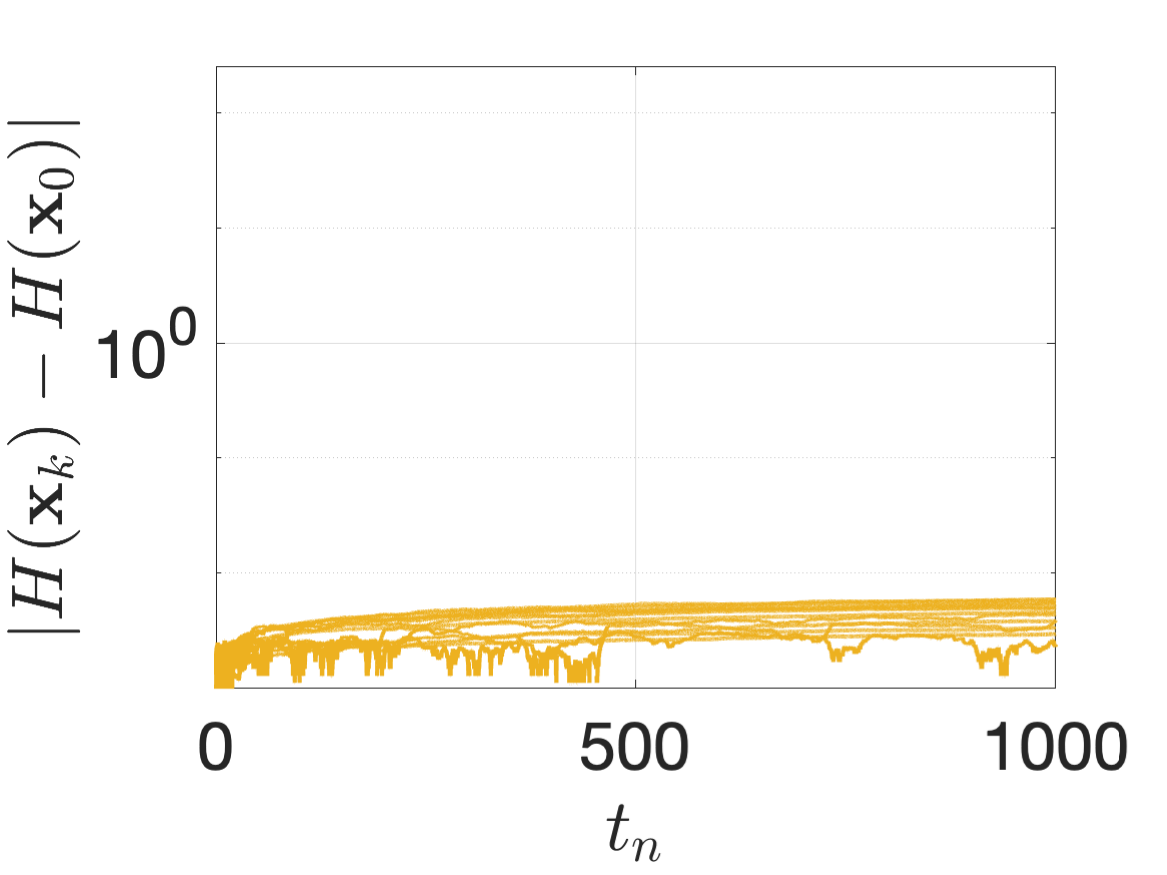}
			\caption{MQAV-MP2}
			\label{fig:quartic_ham_err3}
		\end{subfigure}
		\caption{Numerical orbits of the quartic Hamiltonian ODE \eqref{qodeeg} (first row) and the corresponding errors to the Hamiltonian (second row). The columns corresponds to the midpoint rule (first column), Simpson's rule (second column) and the MQAV midpoint method (third column). The black and red dots are the initial conditions, red meaning that the numerical solution becomes unstable during the $n=10^4$ time steps. The time step is $h=1/10$.}
		\label{fig:quartic_ham}
	\end{figure}
	Numerical simulations using the MQAV midpoint method (MQAV-MP2) are presented in figure \ref{fig:quartic_ham}. As the ODE is Hamiltonian, there exist a number of methods with good structure-preserving properties for comparison. We will use the standard midpoint method (MP2) as well as AVF method. The midpoint method is symplectic and therefore has bounded Hamiltonian error. The AVF method preserves the Hamiltonian and is also a Runge-Kutta method known as Simpson's method \cite{celledoni2009energy} corresponding to evaluation of the integral in \eqref{AVF} with Gaussian quadrature of order four. Furthermore, we remark that the MQAV-MP2 method with $\alpha = 1/3$ also yields the AVF method. All methods are implemented using fixed point iterations with an absolute tolerance of $1.11\times 10^{-15}$. Using a time step of $h=1/10$ a total of 13 numerical orbits are plotted from the initial conditions $\bx_0^{(i)} = (2 + 2i/3,0)^T$, where $i=0,...,12$, denoted by black dots. Red dots denote those initial conditions where the solution becomes unstable. We see that the midpoint method loses stability for the initial conditions corresponding to $i=8,...,12$, the AVF method loses stability for $i=10,...,12$ while the MQAV midpoint method retains stability for all initial conditions. We note that all three methods are unstable when $i=13$. Figure \ref{fig:quartic_ham} also shows the Hamiltonian error as a function of simulation time. We see that the MQAV midpoint method and the Simpson method preserves the Hamiltonian up to machine precision for their stable orbits, whereas the standard midpoint rule does not. Lastly, we note that different choices of $\alpha$ affect the stability and accuracy of the method. 
	
	\subsection{A Nambu system with two integrals}
	In this next example, we have designed an ODE with bounded orbits possessing two integrals, one of degree four and one of degree eight. The integrals are
	\begin{align}
		H_1(\bx) =& {x_{1}}^4\,{x_{2}}^4+x_{1}\,x_{3}+{x_{2}}^4\,{x_{3}}^2 \\
		H_2(\bx) =&\left({x_{2}}^2-1\right)\,\left({x_{1}}^2+{x_{2}}^2+{x_{3}}^2\right)
	\end{align} 
	and the ODE is given by 
	\begin{equation}\label{2intODE}
		\dot{\bx} = \nabla H_1(\bx) \times \nabla H_2(\bx)
	\end{equation}
%
	Now consider the following reduced-degree integrals 
	\begin{align}
		\tilde{H}_1(\bx,\by^{[1]},\by^{[2]}) =& {y^{[2]}_{1,1,1,1}}\,{y^{[2]}_{2,2,2,2}}+x_{1}x_3+{y^{[2]}_{2,2,2,2}}\,{y^{[1]}_{3,3}} \\
		\tilde{H}_2(\bx,\by^{[1]},\by^{[2]}) =&\left({y^{[1]}_{2,2}}-1\right)\,\left({y^{[1]}_{1,1}}+{y^{[1]}_{2,2}}+{y^{[1]}_{3,3}}\right)
	\end{align}
	where $\by^{[1]} = (y^{[1]}_{1,1},y^{[1]}_{2,2},y^{[1]}_{3,3})^T$ and $\by^{[2]} = (y^{[2]}_{1,1,1,1},y^{[2]}_{2,2,2,2},y^{[2]}_{3,3,3,3})^T$. This yields the following MQAV vector field of $\bff(\bx)$
	\begin{align}
		\dot{\bx} =& \nabla \tilde{H}_1(\bx,\by^{[1]},\by^{[2]}) \times \nabla \tilde{H}_2(\bx,\by^{[1]},\by^{[2]})\\ =&\left(\begin{array}{c} 8\,x_2\,x_3\,y^{[1]}_{2,2}\,(y^{[1]}_{3,3} + y^{[2]}_{1,1,1,1})\,(y^{[1]}_{2,2} - 1) - 2\,x_2\,(x_1 + 2\,x_3\,y^{[2]}_{2,2,2,2})\,(y^{[1]}_{1,1} + 2\,y^{[1]}_{2,2} + y^{[1]}_{3,3} - 1)\\
			2\,x_1\,(x_1 + 2\,x_3\,y^{[2]}_{2,2,2,2})\,(y^{[1]}_{2,2} - 1) - 2\,x_3\,(x_3 + 4\,x_1\,y^{[1]}_{1,1}\,y^{[2]}_{2,2,2,2})\,(y^{[1]}_{2,2} - 1)\\
			2\,x_2\,(x_3 + 4\,x_1\,y^{[1]}_{1,1}\,y^{[2]}_{2,2,2,2})\,(y^{[1]}_{1,1} + 2\,y^{[1]}_{2,2} + y^{[1]}_{3,3} - 1) - 8\,x_1\,x_2\,y^{[1]}_{2,2}\,(y^{[1]}_{3,3} + y^{[2]}_{1,1,1,1})\,(y^{[1]}_{2,2} - 1) \end{array}\right).
	\end{align}
	The MQAV midpoint method is now applied, which does not coincide with the AVF method in this case for this choice of $\tilde{H}_1$ and $\tilde{H}_1$. In addition, we also consider the following higher-order methods: a three-stage fourth-order DISRK method (DISRK4), an eleven-stage sixth-order DISRK method (DISRK6) \cite{kalogiratou2014sixth} and a 15-stage eighth-order composition method (C8) \cite{suzuki1993higher}. The values for the substeps $b_i$ are presented in appendix \ref{sec:coeffs}. \\
	
\begin{figure}[h!]
		\centering
		\begin{subfigure}{0.45\textwidth}
			\includegraphics[width=\linewidth]{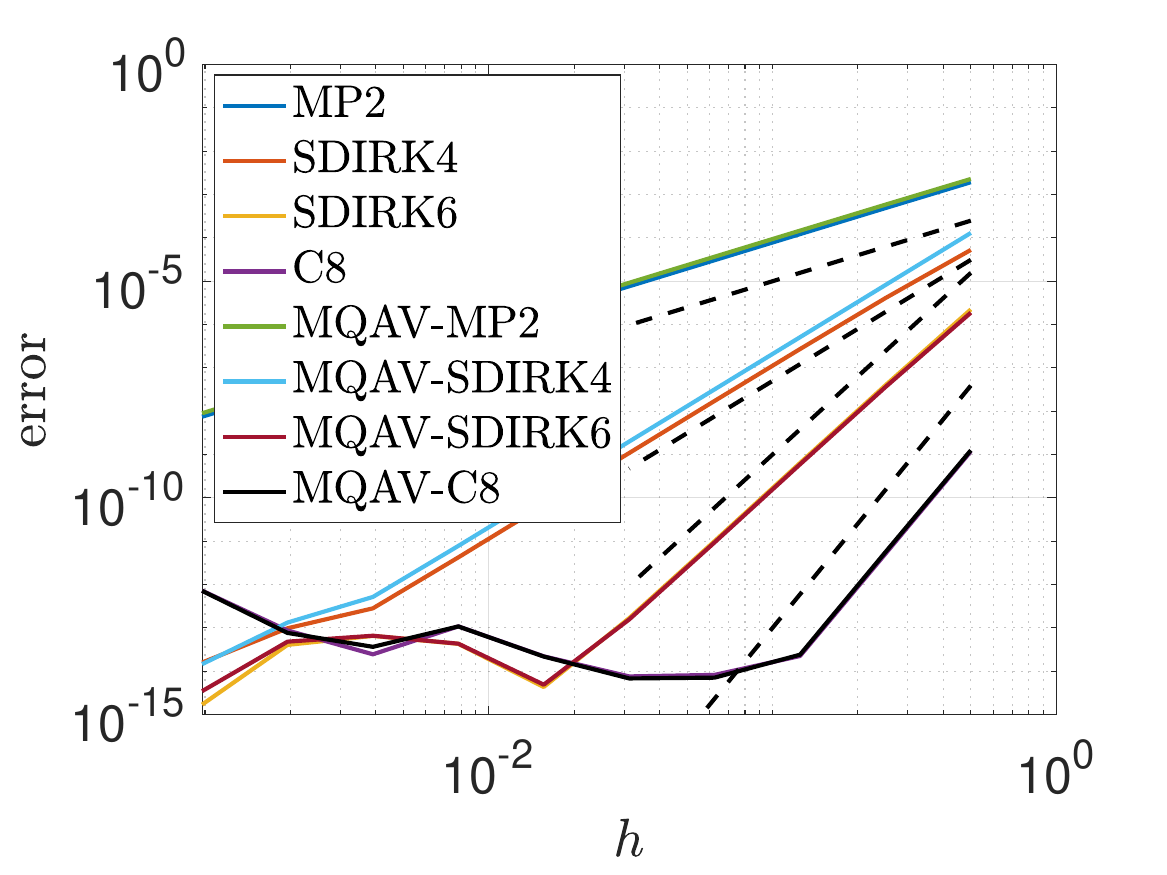}
			\caption{}
			\label{fig:convergence}
		\end{subfigure}
		\begin{subfigure}{0.45\textwidth}
			\includegraphics[width=\linewidth]{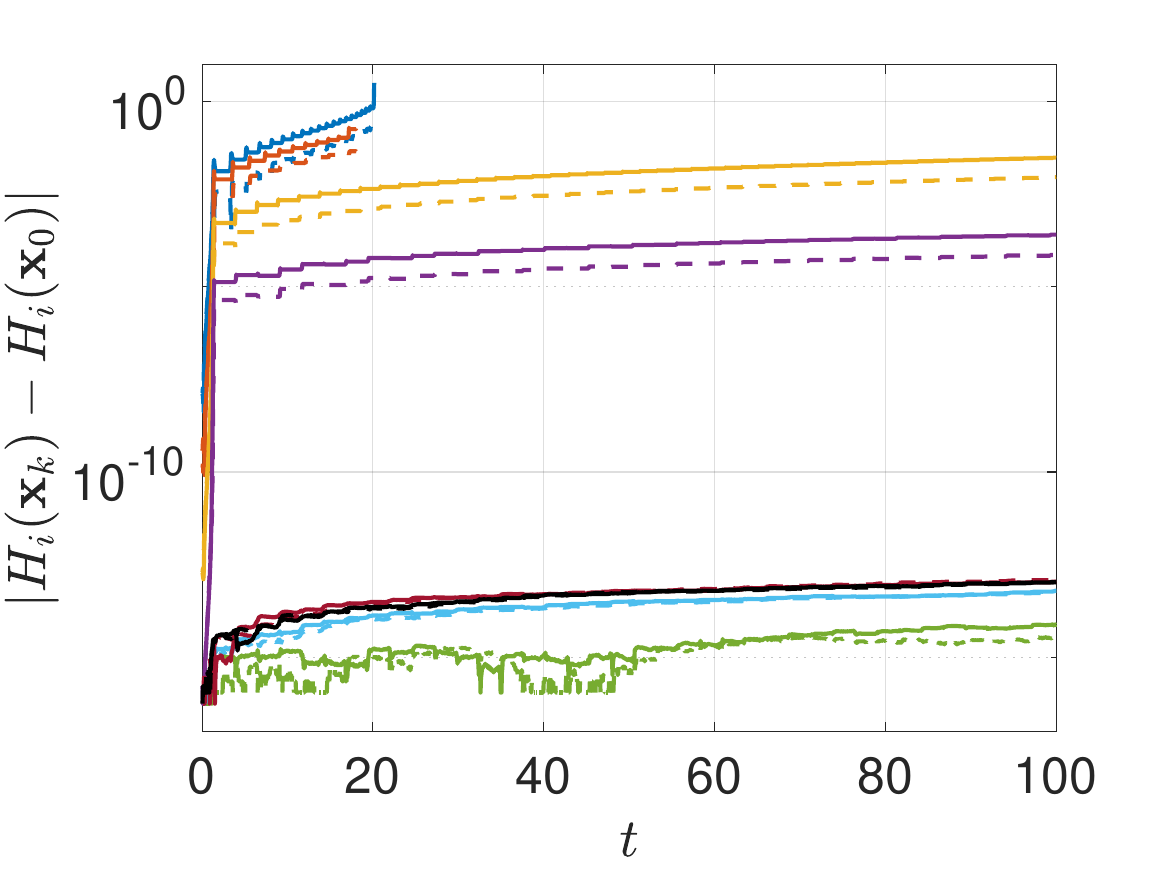}
			\caption{}
			\label{fig:H1_error}
		\end{subfigure}
		\caption{The convergence of the eight methods after one second of simulation time for varying $h$ (a) and the errors $H_1(\bx_k)-H_1(\bx_0)$ (solid lines) and $H_2(\bx_k)-H_2(\bx_0)$ (dashed lines) during 100 seconds of simulation time for $h=1/20$ (b). The black dashed lines in figure (a) are orders two, four, six and eight.}\label{}
	\end{figure}

\begin{figure}[h!]
		\centering
		\begin{subfigure}{0.32\textwidth}
			\includegraphics[width=\linewidth]{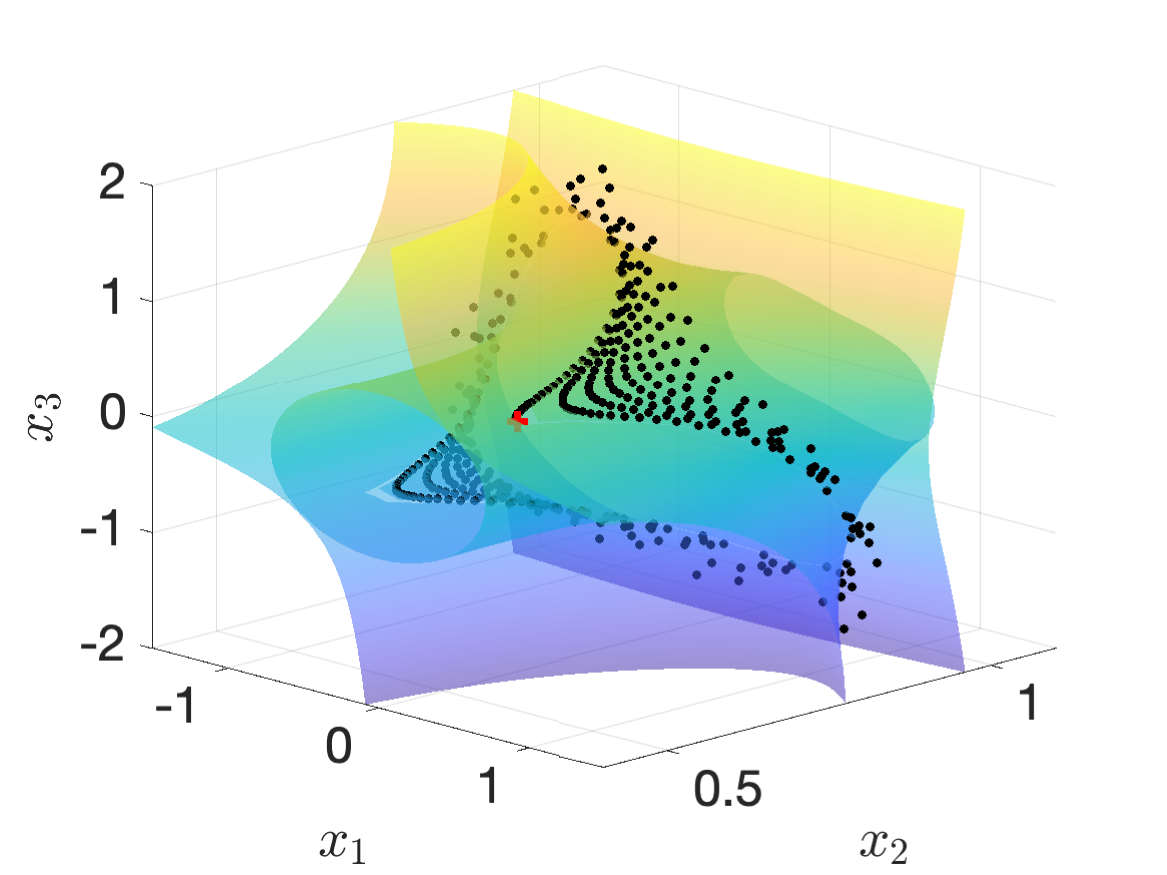}
			\caption{MP2}
			\label{fig:phase_lines_1}
		\end{subfigure}
		\begin{subfigure}{0.32\textwidth}
			\includegraphics[width=\linewidth]{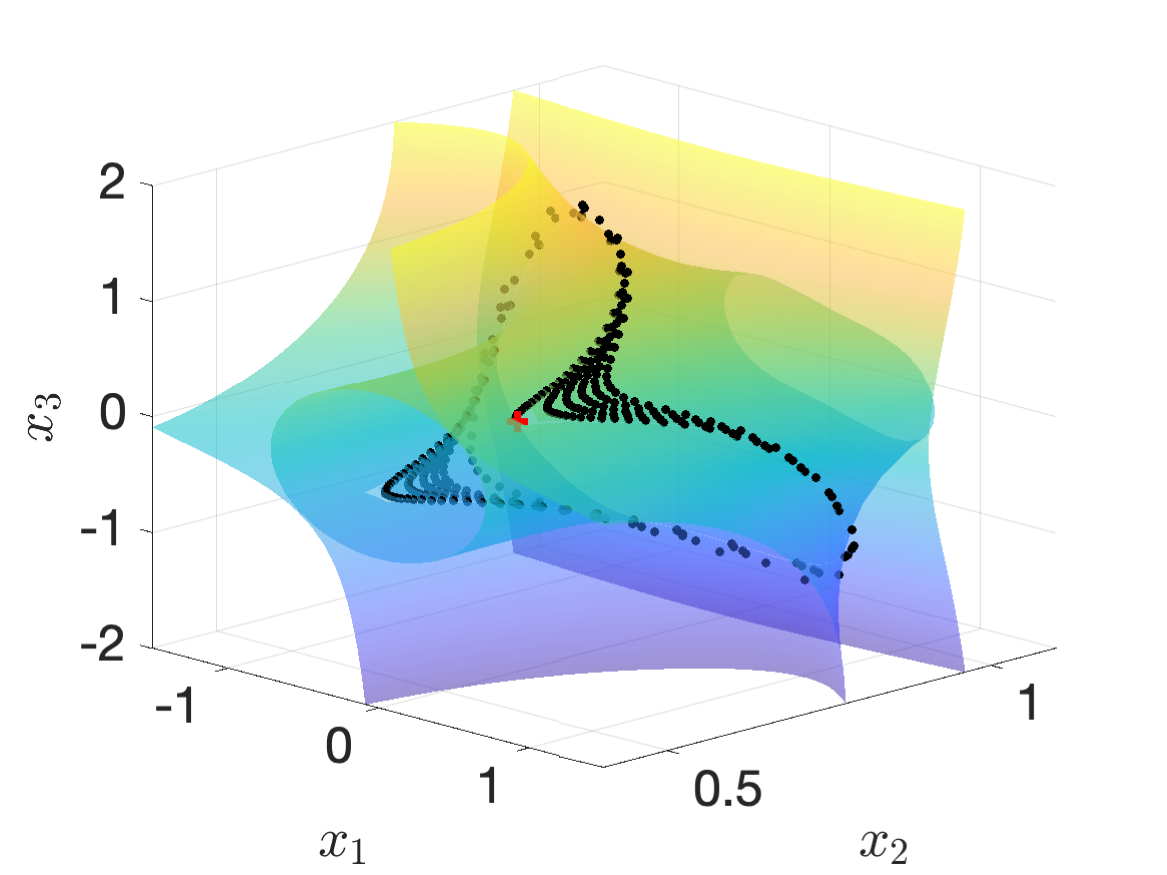}
			\caption{DISRK4}
			\label{fig:phase_lines_2}
		\end{subfigure}
		\begin{subfigure}{0.32\textwidth}
			\includegraphics[width=\linewidth]{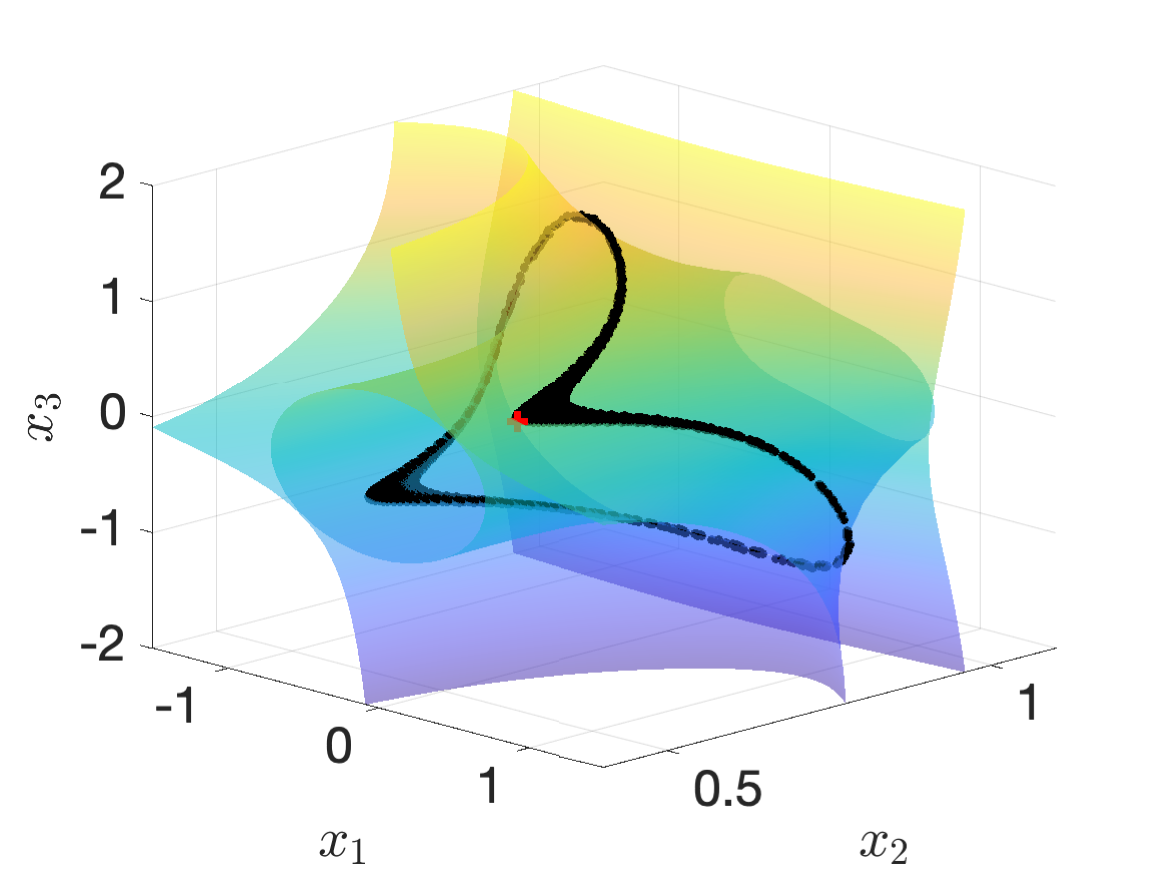}
			\caption{DISRK6}
			\label{fig:phase_lines_3}
		\end{subfigure}
		
		\begin{subfigure}{0.32\textwidth}
			\includegraphics[width=\linewidth]{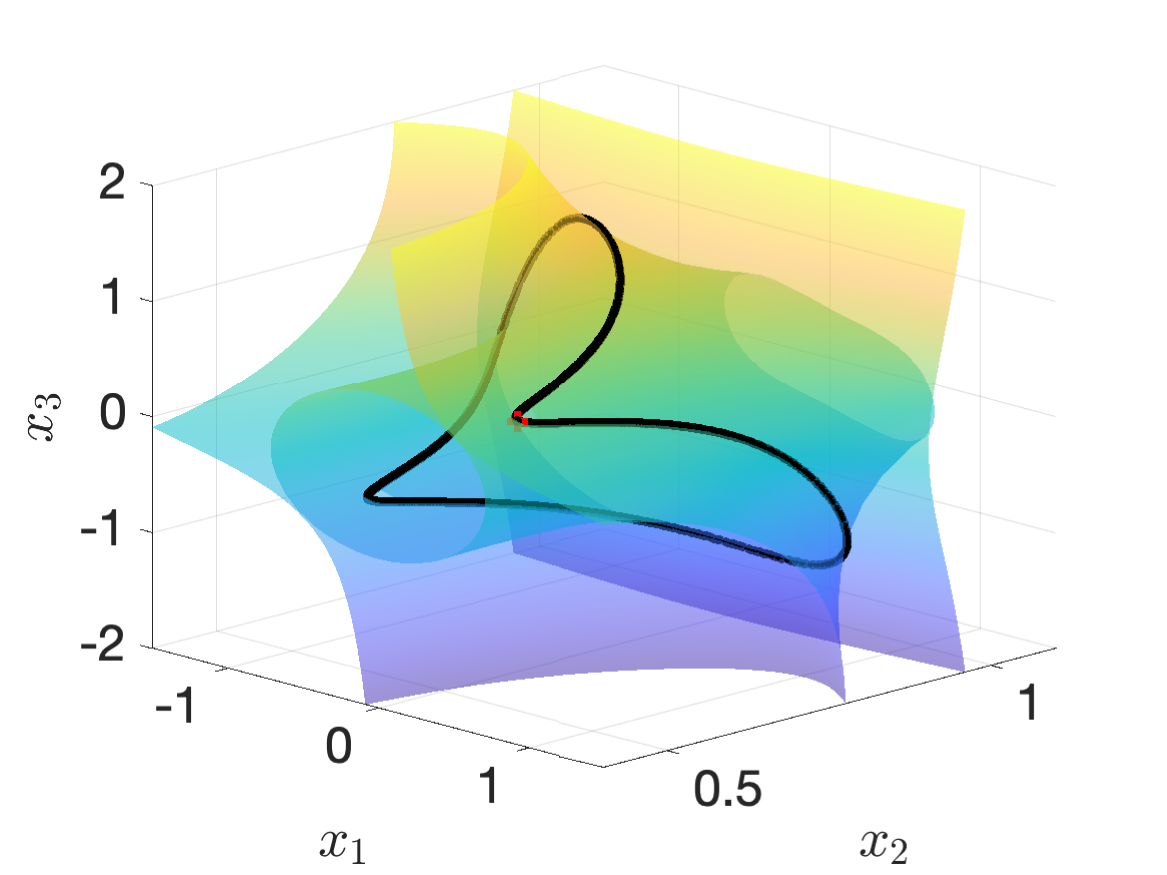}
			\caption{MQAV-MP2}
			\label{fig:phase_lines_4}
		\end{subfigure}
		\begin{subfigure}{0.32\textwidth}
			\includegraphics[width=\linewidth]{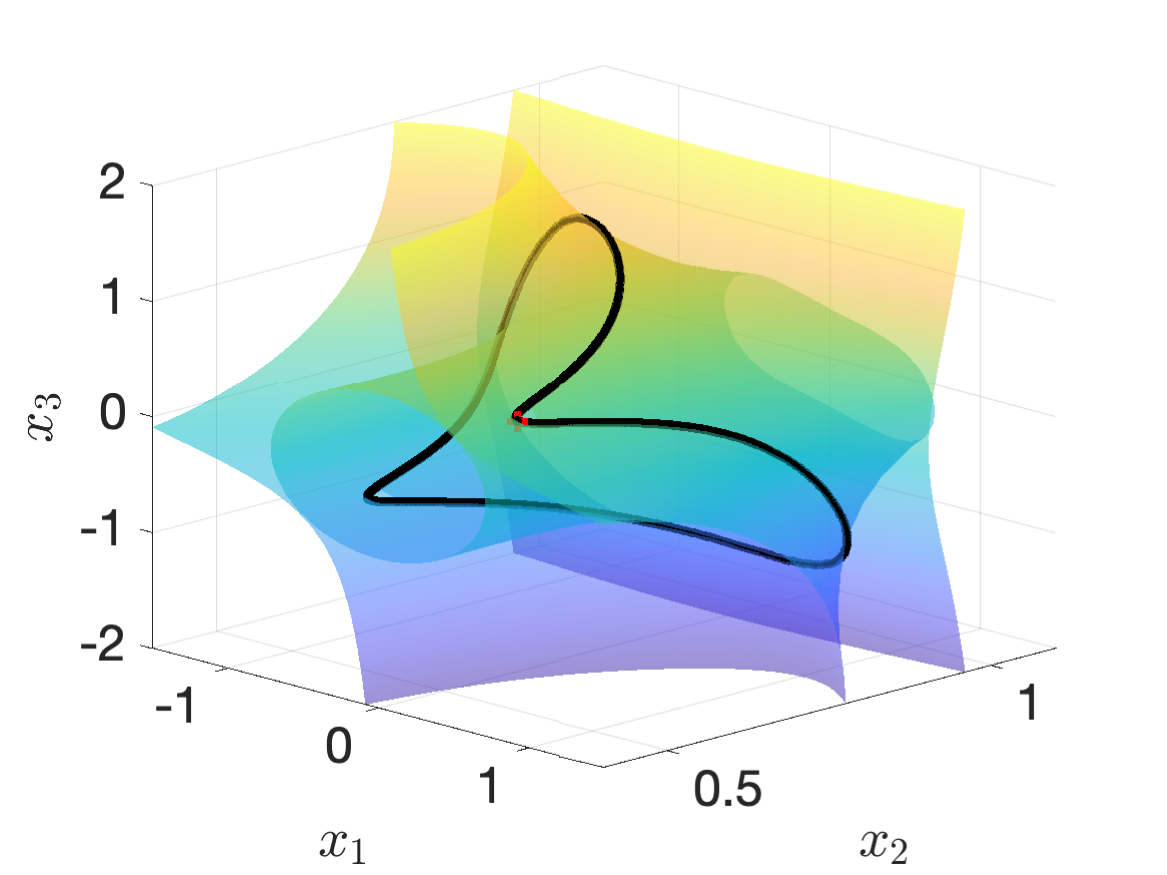}
			\caption{MQAV-DISRK4}
			\label{fig:phase_lines_5}
		\end{subfigure}
		\begin{subfigure}{0.32\textwidth}
			\includegraphics[width=\linewidth]{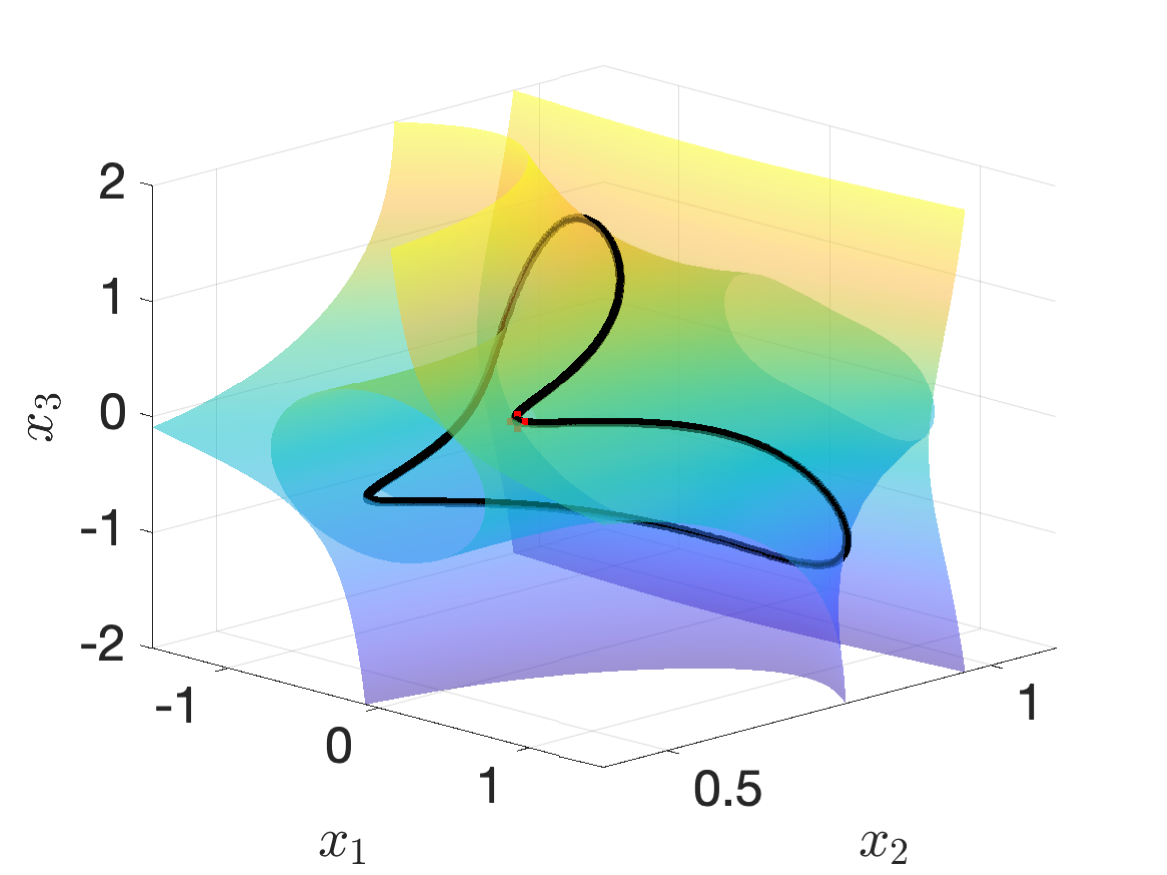}
			\caption{MQAV-DISRK6}
			\label{fig:phase_lines_6}
		\end{subfigure}
		\caption{Numerical phase lines of the ODE using the numerical methods starting from $\bx_0 = (1/2,1/2,1/2)^T$ with $h=1/20$. The composition methods (C8 and MQAV-C8) yield visually similar figures to that figures (d), (e) and (f). }
		\label{fig:phase lines}
	\end{figure}

	A total of eight numerical methods are implemented. The first four are the conventional MP2, DISRK4, DISRK6 and C8 methods applied to the original ODE \eqref{2intODE}. The latter four methods are the MQAV versions of these methods, denoted by MQAV-MP2, MQAV-DISRK4, MQAV-DISRK6 and MQAV-C8. Using the initial condition $\bx_0=(1/2,1/2,1/2)^T$ over the time interval $[0,1]$ the convergence of the six methods are tested and the results are displayed in figure \ref{fig:convergence}. Here, we observe the expected orders. The numerical solutions are calculated from the same initial conditions over the time interval $[0,100]$ with a stepsize of $h = 1/20$. The numerical solutions are plotted in figure \ref{fig:phase lines}. Here we see that the MQAV methods remain on the intersection of the two isosurfaces $H_1(\bx)=H_1(\bx_0)$ and $H_2(\bx)=H_2(\bx_0)$, whereas the conventional methods drift off. The integral errors $H_i(\bx_k)-H_i(\bx_0)$ are presented in figure \ref{fig:H1_error}. We see that the MQAV methods preserve the integrals within machine precision, while the conventional methods do not. Furthermore, the MP2 and DISRK4 methods lose stability early in the simulation.  
	
	%
	
	\subsection{Toda Lattice ($N=3$)}
\begin{figure}[h!]
		\centering
		\includegraphics[width=0.45\linewidth]{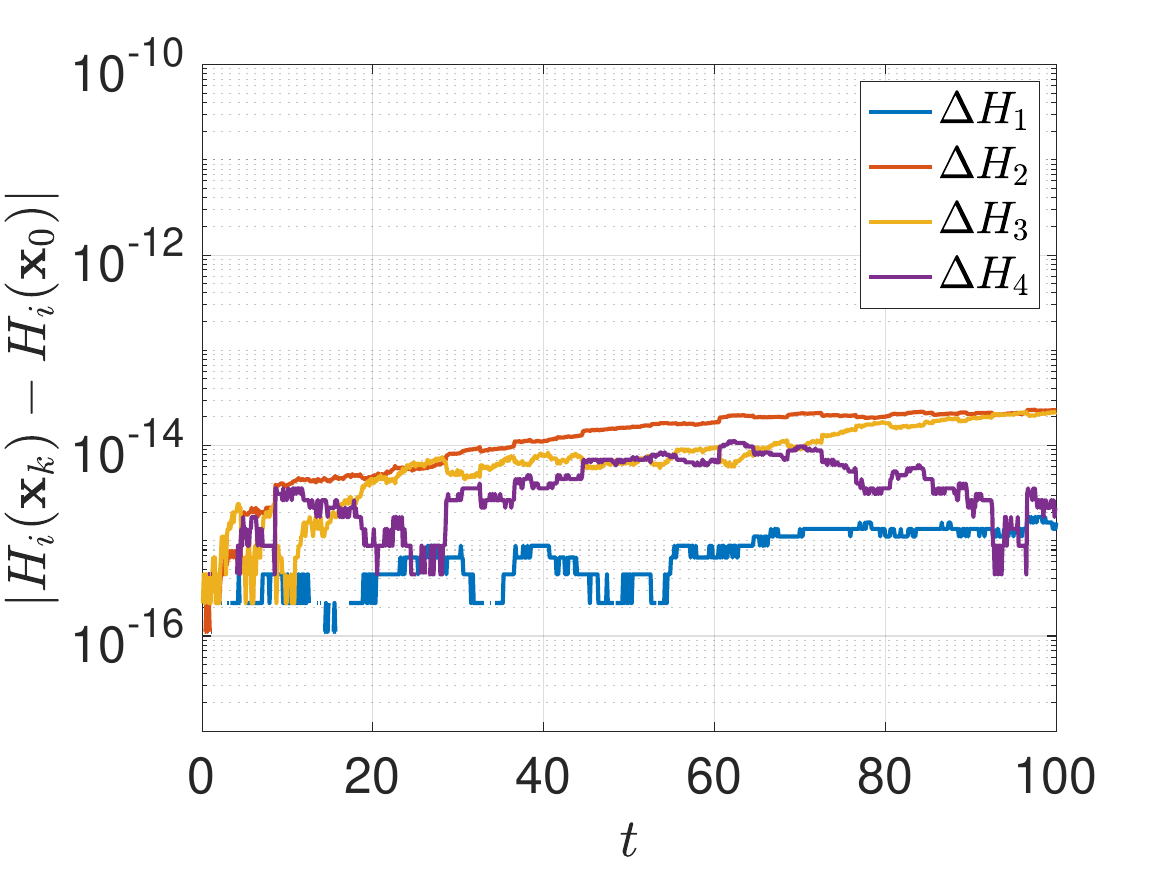}
		\caption{The errors of the integrals of the $N=3$ Toda lattice. }
		\label{fig:todaherrors}
	\end{figure}

	One of the most famous integrable systems in physics is the Toda lattice. Integrable disretisations of this system have been studied in, for example \cite{suris2012problem}. Consider a Toda lattice of $N=3$ particles with periodic boundary conditions in the variables $\bx = (a_1,a_2,a_3,b_1,b_2,b_3)^T$. This is given by \cite{suris2012problem}
	\begin{align}
		\dot{a}_1 =& a_{{1}} \left( b_{{2}}-b_{{1}} \right),  \\
		\dot{a}_2 =& a_{{2}} \left( b_{{3}}-b_{{2}} \right),  \\
		\dot{a}_3 =& a_{{3}} \left( b_{{1}}-b_{{3}} \right),  \\
		\dot{b}_1 =&  a_{{1}}-a_{{3}} , \\
		\dot{b}_2 =&  a_{{2}}-a_{{1}},  \\
		\dot{b}_3 =&  a_{{3}}-a_{{2}}.
	\end{align}
	This ODE possesses the following integrals 
	\begin{align}
		H_1 =& b_{{1}}+b_{{2}}+b_{{3}},\\
		H_2 =& a_1a_2a_3,\\
		H_3 =& \frac{1}{3}\left({ {{b_{{3}}}^{3}}}+{{{b_{{1}}}^{3}}}+{ {{b_
					{{2}}}^{3}}}\right)+a_{{1}}b_{{1}}+a_{{2}}b_{{2}}+a_{{3}}b_{{3}}+a_{{1}}b_{{2}}+a_{{2}}b_{{3}}+a_{{3}}b_{{1}},\\
		H_4 =&  \frac{1}{2}\left({{b_{{1}}}^{2}}+{b_{{2}}}^{2}+
		{{b_{{3}}}^{2}}\right) + a_{{1}}+a_{{2}}+a_{{3}} ,
	\end{align}
	two of which are cubic. Quadratic reduced-degree integrals can be written by introducing the variables $y_{1,2} = a_1a_2$, $y_{4,4} = b_1^2$, $y_{5,5} = b_2^2$ and $y_{6,6} = b_3^2$. This particular choice makes the MQAV method distinct from the AVF method. The skew-symmetric 5-tensor components $S_{i_0,i_1,i_2,i_3,i_4}$ are found using equation \eqref{s}. Using a time step of $h=1/10$ starting from initial conditions $\bx_0 = (1,2,3,4,5,6)/6$ we implement the MQAV midpoint method. Figure \ref{fig:todaherrors} shows the errors of the integrals $\Delta H_i = |H_i(\bx_k)-H_i(\bx_0)|$ for $i=1,...,4$. We see that all four integrals are preserved to machine precision. 
	
	\section{Discussion and conclusion}\label{conclusion}
	In this paper, we have introduced the MQAV method, which is a novel discrete gradient method for preserving polynomial first integrals. The MQAV midpoint method is a Runge-Kutta method but applied to a higher dimensional ODE, whose solution solves the original ODE under study. Due to this, the method inherits the same properties as the underlying Runge-Kutta method on the extended space. In the case of the midpoint rule such properties include time-symmetry, affine equivariance, A-stability and a B-series expansion; however, what can be investigated further is how these properties on the extended phase space relate to the original ODE. For Hamiltonian ODEs, we have observed that the MQAV approach yields the AVF method as a special case. As mentioned, a limitation of the method is when applied to high-dimensional systems with many first integrals. In this setting the computational cost of the method is dominated by the calculation of the $k+1$-skew symmetric tensor $S$. This restricts the applicability of the MQAV method to semi-discretised PDEs with many first integrals, although such problems are not common throughout the literature. In any case, it is desirable to implement an efficient calculation of this tensor in practise. \\
	
	While the primary goal of the paper has been solving conservative ODE systems, we note that the MQAV method can be straightforwardly applied to Hamiltonian PDEs (i.e., the $k=1$ case) for function $u:\R^n\times\R\rightarrow\R^m$
	\begin{equation}
		\dot{u}=\mathcal{S}\frac{\delta \mathcal{H}}{\delta u}, \quad \mathcal{H}[u] = \int \bar{\mathcal{H}}(x;u_x,u_{xx},...) \dd x
	\end{equation}
	where $\mathcal{S}$ is a constant, skew-symmetric, linear differential operator and $\frac{\delta}{\delta u}$ the variational derivative for the Hamiltonian $\mathcal{H}$. See, for example, \cite{celledoni2012preserving} for details. In parallel developments, similar methods also based on the idea of reducing the degree of a cubic Hamiltonian to a quadratic has been implemented to develop high order methods for the quantum Zakharov system \cite{zhang2022arbitrary} and KdV equation \cite{gong2021new}. Furthermore, the authors of \cite{celledoni2012preserving} present a systematic approach to discretising Hamiltonian PDEs into Hamiltonian ODEs. Here, the idea is to semi-discretise the PDE by approximating the skew-symmetric differential operator with a finite dimensional skew-symmetric matrix $\mathcal{S}\rightarrow S$, and the Hamiltonian  with a semi-discretised version of it $\mathcal{H}\rightarrow H$ such that the resulting semi-discretised ODE is in skew-gradient form. That is, in a form amenable to the MQAV method. These promising examples suggest that the MQAV method, in combination with efficient structure-preserving implementation similar to those proposed in \cite{zhang2022arbitrary} can be used to develop efficient methods for a wide range of conservative Hamiltonian PDEs, a topic that we hope will be studied in the future. \\
	
	We conclude by outlining some possible future avenues of research. The fact that the MQAV method yields the AVF method infers that it is also a Runge-Kutta method also on the original phase space. However, this is not necessarily the case for ODEs with non-constant Poisson structure or for arbitrary choice of the free parameters. An interesting question that stems from this is for what choice of free parameters and for which classes of ODE does the MQAV midpoint rule yield a Runge-Kutta method on the original phase space. In this paper we have formulated the MQAV method based on introducing MQAVs of monomial type, that is, of the form $y_{i,j}=x_ix_j$, whereas the method would work just as well by introducing carefully chosen polynomial MQAVs of the form $y = \sum \alpha_i x_{a_i}x_{a_j}$. Such a formulation could lead to a less implicit method and perhaps improve performance. Furthermore, the idea of introducing auxiliary variables to transform an ODE in to a higher dimensional one with simpler features can be applied to ODEs with geometric properties other than polynomial first integrals. For example in \cite{mclachlan1999geometric} the authors show how to preserve \textit{Lyapunov functions} with discrete gradients, which also fits the framework of the MQAV method. In \cite{tapley2021preservation} it is shown that all Runge-Kutta methods preserve all affine \textit{second} integrals  (also known as Darboux polynomials) and \textit{rational affine first integrals}. Therefore, by introducing auxiliary variables to transform second integrals that are polynomial in $\bx$ into second integrals that are linear in $\bx$ and $\by$, one could develop methods that preserve second integrals. In a similar fashion, one could develop methods that preserves arbitrary rational integrals, though it is not clear if such methods are also discrete gradients. 
	
	\section*{Acknowledgments}
	The author would like to thank Brynjulf Owren for his helpful comments and insightful discussions, including pointing out a connection to the AVF method. This work was funded by the European Union’s Horizon 2020 research and innovation program under the Marie Skłodowska-Curie Grant agreement No. 691070. 
	
	\newpage
	\bibliographystyle{unsrt}
	\bibliography{bibliography.bib}

\begin{thebibliography}{10}

\bibitem{mclachlan1998unified}
Robert~I McLachlan, GRW Quispel, and Nicolas Robidoux.
\newblock Unified approach to {H}amiltonian systems, {P}oisson systems,
  gradient systems, and systems with {L}yapunov functions or first integrals.
\newblock {\em Physical Review Letters}, 81(12):2399, 1998.

\bibitem{hairer2006geometric}
Ernst Hairer, Christian Lubich, and Gerhard Wanner.
\newblock {\em Geometric numerical integration: structure-preserving algorithms
  for ordinary differential equations}, volume~31.
\newblock Springer Science \& Business Media, 2006.

\bibitem{brugnano2015reprint}
Luigi Brugnano, Felice Iavernaro, and Donato Trigiante.
\newblock Reprint of analysis of {H}amiltonian boundary value methods
  ({HBVM}s): A class of energy-preserving {R}unge--{K}utta methods for the
  numerical solution of polynomial {H}amiltonian systems.
\newblock {\em Communications in Nonlinear Science and Numerical Simulation},
  21(1-3):34--51, 2015.

\bibitem{miyatake2014energy}
Yuto Miyatake.
\newblock An energy-preserving exponentially-fitted continuous stage
  {R}unge--{K}utta method for {H}amiltonian systems.
\newblock {\em BIT Numerical Mathematics}, 54(3):777--799, 2014.

\bibitem{ranocha2020relaxation}
Hendrik Ranocha and David~I Ketcheson.
\newblock Relaxation {R}unge--{K}utta methods for {H}amiltonian problems.
\newblock {\em Journal of Scientific Computing}, 84(1):1--27, 2020.

\bibitem{iserles2000preserving}
Arieh Iserles and Antonella Zanna.
\newblock Preserving algebraic invariants with {R}unge--{K}utta methods.
\newblock {\em Journal of computational and applied mathematics},
  125(1-2):69--81, 2000.

\bibitem{gonzalez1996time}
Oscar Gonzalez.
\newblock Time integration and discrete {H}amiltonian systems.
\newblock {\em Journal of Nonlinear Science}, 6(5):449--467, 1996.

\bibitem{mclachlan1999geometric}
Robert~I McLachlan, G~R~W Quispel, and Nicolas Robidoux.
\newblock Geometric integration using discrete gradients.
\newblock {\em Philosophical Transactions of the Royal Society of London.
  Series A: Mathematical, Physical and Engineering Sciences},
  357(1754):1021--1045, 1999.

\bibitem{itoh1988hamiltonian}
Toshiaki Itoh and Kanji Abe.
\newblock {H}amiltonian-conserving discrete canonical equations based on
  variational difference quotients.
\newblock {\em Journal of Computational Physics}, 76(1):85--102, 1988.

\bibitem{quispel2008new}
GRW Quispel and David~Ian McLaren.
\newblock A new class of energy-preserving numerical integration methods.
\newblock {\em Journal of Physics A: Mathematical and Theoretical},
  41(4):045206, 2008.

\bibitem{hairer2010energy}
Ernst Hairer.
\newblock Energy-preserving variant of collocation methods.
\newblock {\em Journal of Numerical Analysis, Industrial and Applied
  Mathematics}, 5:73--84, 2010.

\bibitem{celledoni2009energy}
Elena Celledoni, Robert~I McLachlan, David~I McLaren, Brynjulf Owren, G~R~W
  Quispel, and William~M Wright.
\newblock Energy-preserving {R}unge-{K}utta methods.
\newblock {\em ESAIM: Mathematical Modelling and Numerical Analysis},
  43(4):645--649, 2009.

\bibitem{matsuo2001dissipative}
Takayasu Matsuo and Daisuke Furihata.
\newblock Dissipative or conservative finite-difference schemes for
  complex-valued nonlinear partial differential equations.
\newblock {\em Journal of Computational Physics}, 171(2):425--447, 2001.

\bibitem{furihata2019discrete}
Daisuke Furihata and Takayasu Matsuo.
\newblock {\em Discrete variational derivative method: a structure-preserving
  numerical method for partial differential equations}.
\newblock Chapman and Hall/CRC, 2019.

\bibitem{celledoni2012preserving}
Elena Celledoni, Volker Grimm, Robert~I McLachlan, DI~McLaren, D~O’Neale,
  Brynjulf Owren, and GRW Quispel.
\newblock Preserving energy resp. dissipation in numerical {PDE}s using the
  “{A}verage {V}ector {F}ield” method.
\newblock {\em Journal of Computational Physics}, 231(20):6770--6789, 2012.

\bibitem{gong2014some}
Yuezheng Gong, Jiaxiang Cai, and Yushun Wang.
\newblock Some new structure-preserving algorithms for general multi-symplectic
  formulations of {H}amiltonian {PDEs}.
\newblock {\em Journal of Computational Physics}, 279:80--102, 2014.

\bibitem{li2015general}
Yu-Wen Li and Xinyuan Wu.
\newblock General local energy-preserving integrators for solving
  multi-symplectic {H}amiltonian {PDE}s.
\newblock {\em Journal of Computational Physics}, 301:141--166, 2015.

\bibitem{dahlby2011general}
Morten Dahlby and Brynjulf Owren.
\newblock A general framework for deriving integral preserving numerical
  methods for {PDE}s.
\newblock {\em SIAM Journal on Scientific Computing}, 33(5):2318--2340, 2011.

\bibitem{eidnes2019linearly}
S{\o}lve Eidnes, Lu~Li, and Shun Sato.
\newblock Linearly implicit structure-preserving schemes for {H}amiltonian
  systems.
\newblock {\em Journal of Computational and Applied Mathematics}, page 112489,
  2019.

\bibitem{shen2018scalar}
Jie Shen, Jie Xu, and Jiang Yang.
\newblock The scalar auxiliary variable (sav) approach for gradient flows.
\newblock {\em Journal of Computational Physics}, 353:407--416, 2018.

\bibitem{kemmochi2021scalar}
Tomoya Kemmochi and Shun Sato.
\newblock Scalar auxiliary variable approach for conservative/dissipative
  partial differential equations with unbounded energy.
\newblock {\em arXiv preprint arXiv:2105.04055}, 2021.

\bibitem{yang2017numerical}
Xiaofeng Yang, Jia Zhao, Qi~Wang, and Jie Shen.
\newblock Numerical approximations for a three-component {Cahn--Hilliard}
  phase-field model based on the invariant energy quadratization method.
\newblock {\em Mathematical Models and Methods in Applied Sciences},
  27(11):1993--2030, 2017.

\bibitem{gong2021new}
Yuezheng Gong, Yue Chen, Chuwu Wang, and Qi~Hong.
\newblock A new class of high-order energy-preserving schemes for the
  {Korteweg-de Vries} equation based on the quadratic auxiliary variable
  ({QAV}) approach.
\newblock {\em arXiv e-prints}, pages arXiv--2108, 2021.

\bibitem{cooper1987stability}
GJ~Cooper.
\newblock Stability of {R}unge-{K}utta methods for trajectory problems.
\newblock {\em IMA journal of numerical analysis}, 7(1):1--13, 1987.

\bibitem{celledoni2014minimal}
Elena Celledoni, Brynjulf Owren, and Yajuan Sun.
\newblock The minimal stage, energy preserving {R}unge--{K}utta method for
  polynomial {H}amiltonian systems is the averaged vector field method.
\newblock {\em Mathematics of Computation}, 83(288):1689--1700, 2014.

\bibitem{norton2015projection}
G.~R. W. Quispel Ari~Stern Richard A.~Norton, David I.~McLaren and Antonella
  Zanna.
\newblock Projection methods and discrete gradient methods for preserving first
  integrals of {ODEs}.
\newblock {\em Discrete and Continuous Dynamical Systems}, 35(5):2079--2098,
  2015.

\bibitem{quispel1997solving}
GRW Quispel, HW~Capel, et~al.
\newblock Solving {ODE}’s numerically while preserving all first integrals.
\newblock {\em Preprint, La Trobe University, Melbourne}, 1997.

\bibitem{mclachlan2003spatial}
Robert~I McLachlan.
\newblock Spatial discretization of partial differential equations with
  integrals.
\newblock {\em IMA journal of numerical analysis}, 23(4):645--664, 2003.

\bibitem{mclachlan1995numerical}
Robert~I McLachlan.
\newblock On the numerical integration of ordinary differential equations by
  symmetric composition methods.
\newblock {\em SIAM Journal on Scientific Computing}, 16(1):151--168, 1995.

\bibitem{kalogiratou2014sixth}
Zacharoula Kalogiratou, Theodore Monovasilis, and TE~Simos.
\newblock A sixth order symmetric and symplectic diagonally implicit
  {R}unge-{K}utta method.
\newblock In {\em AIP Conference Proceedings}, volume 1618, pages 833--838.
  American Institute of Physics, 2014.

\bibitem{suzuki1993higher}
M~Suzuki and K~Umeno.
\newblock Higher-order decomposition theory of exponential operators and its
  applications to {QMC} and nonlinear dynamics.
\newblock In {\em Computer simulation studies in condensed-matter physics VI},
  pages 74--86. Springer, 1993.

\bibitem{suris2012problem}
Yuri~B Suris.
\newblock {\em The problem of integrable discretization: {H}amiltonian
  approach}, volume 219.
\newblock Birkh{\"a}user, 2012.

\bibitem{zhang2022arbitrary}
Gengen Zhang and Chaolong Jiang.
\newblock Arbitrary high-order structure-preserving methods for the quantum
  {Z}akharov system.
\newblock {\em arXiv preprint arXiv:2202.13052}, 2022.

\bibitem{tapley2021preservation}
Benjamin~K Tapley.
\newblock On the preservation of second integrals by {R}unge-{K}utta methods.
\newblock {\em arXiv preprint arXiv:2105.10929}, 2021.

\end{thebibliography}


\begin{thebibliography}{10}

\bibitem{hairer2006geometric}
Ernst Hairer, Christian Lubich, and Gerhard Wanner.
\newblock {\em Geometric numerical integration: structure-preserving algorithms
  for ordinary differential equations}, volume~31.
\newblock Springer Science \& Business Media, 2006.

\bibitem{dahlby2011preserving}
Morten Dahlby, Brynjulf Owren, and Takaharu Yaguchi.
\newblock Preserving multiple first integrals by discrete gradients.
\newblock {\em Journal of Physics A: Mathematical and Theoretical},
  44(30):305205, 2011.

\bibitem{calvo2006preservation}
Manuel Calvo, D~Hern{\'a}ndez-Abreu, Juan~I Montijano, and Luis R{\'a}ndez.
\newblock On the preservation of invariants by explicit runge--kutta methods.
\newblock {\em SIAM Journal on Scientific Computing}, 28(3):868--885, 2006.

\bibitem{quispel1997solving}
GRW Quispel, HW~Capel, et~al.
\newblock Solving ode’s numerically while preserving all first integrals.
\newblock {\em Preprint, La Trobe University, Melbourne}, 1997.

\bibitem{mclachlan1999geometric}
Robert~I McLachlan, G~Reinout~W Quispel, and Nicolas Robidoux.
\newblock Geometric integration using discrete gradients.
\newblock {\em Philosophical Transactions of the Royal Society of London.
  Series A: Mathematical, Physical and Engineering Sciences},
  357(1754):1021--1045, 1999.

\bibitem{norton2015projection}
G.~R. W. Quispel Ari Stern Antonella~Zanna Richard A.~Norton, David I.~McLaren.
\newblock Projection methods and discrete gradient methods for preserving first
  integrals of odes.
\newblock {\em Discrete & Continuous Dynamical Systems}, 35(5):2079--2098,
  2015.

\bibitem{mclachlan2002splitting}
Robert~I McLachlan and G~Reinout~W Quispel.
\newblock Splitting methods.
\newblock {\em Acta Numerica}, 11:341, 2002.

\bibitem{cooper1987stability}
GJ967831 Cooper.
\newblock Stability of runge-kutta methods for trajectory problems.
\newblock {\em IMA journal of numerical analysis}, 7(1):1--13, 1987.

\bibitem{mclachlan1998unified}
Robert~I McLachlan, GRW Quispel, and Nicolas Robidoux.
\newblock Unified approach to hamiltonian systems, poisson systems, gradient
  systems, and systems with lyapunov functions or first integrals.
\newblock {\em Physical Review Letters}, 81(12):2399, 1998.

\bibitem{mclachlan1995numerical}
Robert~I McLachlan.
\newblock On the numerical integration of ordinary differential equations by
  symmetric composition methods.
\newblock {\em SIAM Journal on Scientific Computing}, 16(1):151--168, 1995.

\bibitem{celledoni2009energy}
Elena Celledoni, Robert~I McLachlan, David~I McLaren, Brynjulf Owren,
  G~Reinout~W Quispel, and William~M Wright.
\newblock Energy-preserving runge-kutta methods.
\newblock {\em ESAIM: Mathematical Modelling and Numerical Analysis},
  43(4):645--649, 2009.

\bibitem{kalogiratou2014sixth}
Zacharoula Kalogiratou, Theodore Monovasilis, and TE~Simos.
\newblock A sixth order symmetric and symplectic diagonally implicit
  runge-kutta method.
\newblock In {\em AIP Conference Proceedings}, volume 1618, pages 833--838.
  American Institute of Physics, 2014.

\bibitem{suzuki1993higher}
M~Suzuki and K~Umeno.
\newblock Higher-order decomposition theory of exponential operators and its
  applications to qmc and nonlinear dynamics.
\newblock In {\em Computer simulation studies in condensed-matter physics VI},
  pages 74--86. Springer, 1993.

\bibitem{suris2012problem}
Yuri~B Suris.
\newblock {\em The problem of integrable discretization: Hamiltonian approach},
  volume 219.
\newblock Birkh{\"a}user, 2012.

\bibitem{celledoni2019geometric}
Elena Celledoni, DI~McLaren, Brynjulf Owren, and GRW Quispel.
\newblock Geometric and integrability properties of kahan’s method: the
  preservation of certain quadratic integrals.
\newblock {\em Journal of Physics A: Mathematical and Theoretical},
  52(6):065201, 2019.

\end{thebibliography}
	\newpage
	\appendix

	\end{document}